\documentclass{amsart}
\usepackage{amsmath,bm}
\usepackage{amssymb}
\usepackage{graphicx}
\usepackage{mathtools}
\usepackage{tabularx}     
\usepackage{siunitx}  
\usepackage{booktabs}     
\usepackage[pagewise]{lineno}
\usepackage{array}        
\usepackage{caption}      
\usepackage{amsthm}
\newtheoremstyle{boldremark}%
  {3pt}{3pt}%
  {\normalfont}
  {}
  {\bfseries}
  {.}
  {.5em}
  {}

\theoremstyle{boldremark}

\usepackage{float} 
\usepackage{hyperref}
\usepackage{tikz}
\usetikzlibrary{spy}

\usepackage{subcaption}
 
\usepackage{url} 

\usepackage{booktabs}
\title[Data-Driven Filter Design]{Data-Driven Filter Design for Flexible and Noise-Robust Tomographic Imaging}

\author{Hamid Fathi$^{*1,2}$}
\author{Alexander Skorikov$^{1}$}
\author{Tristan van Leeuwen$^{1,2}$}
\date{%
$^{1}$ Computational Imaging Group, Centrum Wiskunde \& Informatica (CWI), Science Park 123, 1098 XG Amsterdam, The Netherlands\\
$^{2}$ Mathematical Institute, Utrecht University, Budapestlaan 6, 3584 CD Utrecht, The Netherlands\\
$^{*}$ Corresponding author.
}
\begin{document}
\begin{abstract}
While filtered back projection (FBP) is still the method of choice for fast tomographic reconstruction, its performance degrades noticeably in the presence of noise, incomplete sampling, or non-standard scan geometries. We propose a data-driven approach for learning FBP filters and projection weights directly from training data, with the goal of improving robustness without sacrificing computational efficiency. The resulting reconstructions adapt naturally to the noise level and acquisition geometry, while retaining the speed and simplicity of classical back-projection. The proposed method can be formulated as a regularized optimization problem for a linear inverse operator, which allows us to establish existence, uniqueness, and stability of the learned solution. From a spectral viewpoint, the learned filters act as data-adaptive gain functions that explicitly balance noise amplification and bias, in close analogy to a regularized pseudo-inverse. Experiments in both 2D and 3D show consistent improvements over conventional FBP and FDK in different case studies. Finally, we show that filters trained on synthetic laminography data generalize well to real-world measurements, delivering image quality comparable to advanced iterative methods without the high computational cost.
\end{abstract}

\maketitle
\section{Introduction}
X-ray computed tomography (CT) is widely used in industrial applications for non-destructive evaluation of objects. The goal is to reconstruct two- or three-dimensional images from X-ray measurements acquired at multiple viewing angles \cite{Kak,Hensen,Gabor}. Well-known examples include medical CT and cone-beam CT. Another important example is X-ray laminography, a technique particularly well suited for imaging flat or layered specimens such as electronic components, where full rotation is often infeasible due to size or mechanical constraints. 

Image reconstruction methods for CT generally fall into two categories: \emph{algebraic approaches}, which iteratively solve a system of linear equations, and \emph{analytical approaches}, which provide a direct solution through filtering and backprojection. 

In the algebraic approach, X-ray measurements are represented through a system of linear equations. Various iterative methods have been developed to solve the linear system \cite{Kak}. A key advantage of these techniques is their flexibility in incorporating prior information or regularization terms to suppress artifacts caused by noise, incomplete data, or other acquisition imperfections \cite{kudo,bilgic}. However, they are often computationally intensive, requiring significant time and memory resources for high-resolution reconstructions \cite{geyer}. Commonly employed iterative algorithms in algebraic approaches include the Algebraic Reconstruction Technique (ART) \cite{ART}, the Simultaneous Iterative Reconstruction Technique (SIRT) \cite{SIRT}, Nesterov accelerated gradient descent (NAG) \cite{Nesterov} and more advanced regularized iterative methods such as Total Variation (TV)-regularized reconstruction \cite{kudo}, which aim to improve image quality and robustness at the expense of additional computation. 

In contrast, analytical approaches are attractive due to their simplicity and computational efficiency \cite{Pan}. They provide a fast, non-iterative solution to the reconstruction problem \cite{Hensen}. A widely used example is the Feldkamp–Davis–Kress (FDK) algorithm \cite{FDK}, which is highly effective in conventional cone-beam CT configurations, particularly when the measurements are of high quality. However, its performance can degrade significantly in the presence of noise \cite{Sidky_2008}. In FDK, standard filters such as Shepp–Logan, Hann, and Hamming use fixed window functions to attenuate high frequencies and suppress artifacts, however they are not tailored to the specific noise characteristics of the measured data \cite{Buzug}. Furthermore, the classical FDK algorithm is derived under the assumption of a circular source trajectory, which limits its direct applicability to more general scanning geometries. Although analytical derivations for laminography and related geometries have been presented \cite{Lauritsch,Myagotin,Yang}, they are mathematically involved, highly dependent on specific acquisition setups, and often lead to strong artifacts due to the missing-cone issue. 

These significant challenges with fixed, pre-defined filters have motivated a different line of research: optimizing the filter itself. This approach seeks to combine the computational efficiency of the filtered back projection (FBP) framework with the superior image quality of iterative methods. The core idea is to move beyond fixed filters (like Shepp-Logan) and instead design filters that are better adapted to the specific geometry and noise statistics.

This concept of an optimized FBP is well-supported in the literature. A theoretical bridge between algebraic and analytical methods was presented in \cite{Older}, and it was shown that any linear shift-invariant algorithm can be represented as an FBP method with some associated filter \cite{Clack_1992}. This led to methods that explicitly derive an FBP filter from an iterative algorithm, allowing FBP to approximate the results of iterative reconstruction \cite{Zeng,Joost,Dan_2}. Others have designed geometry-specific filters for each frequency component in tomosynthesis \cite{Nielsen_2012} or proposed data-dependent filtering strategies \cite{Dan} to reconstruct the image in limited data or noisy scenarios better than classical FBP method. Finally, some works \cite{blumensath} proposed to apply a convolution kernel derived from a specific geometry to deblur the reconstructed image.


We propose a data-driven framework to optimize the filter and (when needed) projection weights in filtered back-projection, with a particular focus on 3D cone-beam circular laminography, where limited-angle sampling makes reconstruction especially challenging. The learning problem is posed as a regularized optimization of a linear inverse operator, which enables existence, uniqueness, and stability guarantees for the learned solution and yields an interpretable spectral-gain view analogous to a regularized pseudo-inverse. We instantiate the framework with geometry-specific parameterizations (parallel-beam, fan-beam, elliptical trajectories, and 3D circular laminography) and demonstrate consistent improvements over classical FBP/FDK, including validation on experimental laminography measurements (LEGO dataset) showing generalization beyond the synthetic training setting.

The outline of the paper is as follows. In Section 2, we introduce the notation and briefly review the inverse problem formulation and fast approximate inversion via filtered back-projection. In Section 3, we present our learning framework for a regularized linear inverse operator, discuss its theoretical properties (existence, uniqueness, and stability), provide a spectral interpretation of the learned filters as data-adaptive gain functions, and detail geometry-informed parameterizations together with implementation aspects. In Section 4, we evaluate the proposed approach across different case studies, including 2D parallel-beam and fan-beam CT, a 2D elliptical trajectory, and 3D circular laminography on both synthetic and experimental data. In Section 5, we conclude the paper and summarize key findings.

\section{Notation and preliminaries}
\subsection{Notation}
Throughout, scalars and vector elements are denoted by lower-case italic letters, vectors by lower-case boldface letters, and matrices by upper-case boldface letters. The notation $\cdot^\top$ stands for the transpose of either a vector or matrix. By $\|\bm{x}\|=\sqrt{\sum_{i=1}^n x_i^2}$ we denote the $\ell_2$-norm of a vector. The Frobenius norm of a matrix $\mathbf{A} \in \mathbb{R}^{m \times n}$ is denoted by $\|\mathbf{A}\|_F = \sqrt{\sum_{i=1}^m \sum_{j=1}^n A_{ij}^2}$. $\operatorname{diag}(\bm{a})$ represents a diagonal matrix with elements of $\bm{a}$ on the main diagonal. $\mathbb{E}[\cdot]$ denotes the expectation. $\operatorname{Tr}(\mathbf{A})$ denotes the trace of a square matrix $\mathbf{A}$, defined as the sum of its diagonal elements. The symbol $\otimes$ represents the Kronecker product, and $\operatorname{vec}(\mathbf{A})$ is the vectorization operator that stacks the columns of a matrix $\mathbf{A}$ into a single column vector.
\subsection{The inverse problem and regularization}
\label{Inverse problem}
The X-ray tomographic acquisition process can be modeled as, 
\begin{align}
\mathbf{A}\mathbf{x}=\mathbf{y},
\end{align}
where $\mathbf{A}\in\mathbb{R}^{M\times N}$ is the forward operator, $\mathbf{x}\in\mathbb{R}^{N}$ represents the object, and $\mathbf{y}\in\mathbb{R}^{M}$ is the measured data. The structure of $\mathbf{A}$ reflects the geometry and ordering of the data acquisition. In a scan comprising $m$ projection angles, each projection is measured by a detector with $d$ pixels, yielding a total of $M=m \cdot d$ rows in $\mathbf{A}$. Each row corresponds to a single X-ray path, modeled as a line integral through the object $\mathbf{x}$ from the source to a specific detector pixel. Rows can be grouped into blocks of $d$ consecutive rows, each block associated with a fixed projection angle.

The goal is to find an $\mathbf{x}$ that best solves this system. 
This inverse problem is typically ill-posed, as $\mathbf{A}$ is ill-conditioned and often there are fewer measurements than unknowns (i.e., $M\leq N$). A basic approach to stabilize the problem is through Tikhonov regularization, which approximates the solution through a regularized pseudo-inverse
\[
\widehat{\mathbf{x}} = \mathbf{A}^\dagger_\lambda \mathbf{y},
\]
where $\mathbf{A}^\dagger_\lambda$ is expressed in terms of the singular value decomposition of $\mathbf{A}$ as
\[
\mathbf{A}^\dagger_\lambda = \mathbf{V}_r\left(\mathbf{\Sigma}_r + \lambda \mathbf{I}\right)^{-1}\mathbf{U}_r^\top,
\]
where $r$ is rank of $\mathbf{A}$, $\mathbf{U}_r\in\mathbb{R}^{M\times r}$ contains the first $r$ left singular vectors of $\mathbf{A}$, $\mathbf{\Sigma}_r \in \mathbb{R}^{r\times r}$ as a diagonal matrix which contains the $r$ largest singular values of $\mathbf{A}$, and $\mathbf{V}_r\in\mathbb{R}^{N\times r}$ contains the first $r$ right singular vectors of $\mathbf{A}$. Furthermore $\lambda \geq  0$ is the regularization parameter. In practice, it is not feasible to explicitly compute the regularized pseudo-inverse and instead the regularized solution is approximated by iteratively solving a regularized least-squares problem \cite{Saad}. 
\subsection{Fast approximate inversion}
For tomographic imaging, filtered back projection offers a fast alternative to the above-mentioned reconstruction methods. They are based on the fact that for many geometries the inverse (on a continuous level) can be expressed as a back projection of filtered data \cite{palamodov2016reconstruction}. Applying the filter through a fast Fourier transform then yields an efficient approximate inverse:
\[
\mathbf{A}^\dagger \approx \mathbf{A}^\top  \mathbf{F}^{-1}\mathbf{D}\mathbf{F},
\]
where $\mathbf{F}\in \mathbb{C}^{M\times M}$ represents a discrete Fourier transform (which can be efficiently applied to vectors using the FFT) and $\mathbf{D}\in \mathbb{R}^{M\times M}$ is a diagonal matrix containing the filter weights. Note that in case $M\leq N$ and $\mathbf{A}$ has full rank we have $\mathbf{A}^\dagger = \mathbf{A}^\top (\mathbf{A}\mathbf{A}^\top)^{-1}$. Regularization is achieved in these methods by designing a filter that cuts off high-frequency contributions. Aside from hand-crafted filters, learned filters have also been considered in the literature \cite{faucris, Tan_DeepFBP}.

\section{A Learned Inverse for tomographic reconstruction}\label{Learned_inverse}

\subsection{Learning a regularized pseudo-inverse}
In this section, we formalize the process of learning a filter. We start with the following more general problem. Given $\mathbf{A} \in \mathbb{R}^{M \times N}$ with $M \le N$, the aim is to learn the linear inverse mapping parametrized by $\mathbf{B}\in \mathbb{R}^{M\times M}$,
\begin{equation}\label{Learned Inverse_supervised}
    \min_{\mathbf{B}~\in~\mathbb{R}^{M\times M}} \;
    \mathbb{E}
    \|\mathbf{A}^\top \mathbf{B} \mathbf{y} - \mathbf{x}\|_2^2 + \lambda \rho(\mathbf{B}),
\end{equation}
where we take the expectation over training pairs $(\mathbf{x},\mathbf{y})$,$\lambda \geq 0$, and $\rho : \mathbb{R}^{M\times M}\rightarrow \mathbb{R}$ is a regularization function. The training pairs typically consist of ground truth objects $\mathbf{x}$ and corresponding noisy measurements $\mathbf{y} \approx \mathbf{A}\mathbf{x}$.

\subsubsection{General Assumptions}To ensure the well-posedness of the problem above, we make the following standing assumptions throughout the analysis:\begin{itemize}\item[\textbf{(A1)}] We assume that the underlying distribution of the training data has finite second moments, i.e., $\mathbb{E}\|\mathbf{x}\|_2^2 < \infty,~ \mathbb{E}\|\mathbf{y}\|_2^2 < \infty$, and the second-order moment matrices $\mathbf{\Sigma}_{\mathbf{xx}} = \mathbb{E}\mathbf{x}\mathbf{x}^\top$ and $\mathbf{\Sigma}_{\mathbf{yy}}=\mathbb{E}\mathbf{y}\mathbf{y}^\top$ are positive definite.
\item[\textbf{(A2)}] The matrix $\mathbf{A}$ has full row rank, implying $\mathbf{A}\mathbf{A}^\top \succ 0$.\item[\textbf{(A3)}] The function $\rho$ is $\mu$-strongly convex.\end{itemize}
With these assumptions in place, we can address the well-posedness of (\ref{Learned Inverse_supervised}).

\subsubsection{Existence, uniqueness, and stability}\label{wellposedness}
Under Assumptions (A1)-(A3), the objective function in \eqref{Learned Inverse_supervised} is strongly convex with respect to $\mathbf{B}$ and hence has a unique solution. Moreover, the optimal solution $\mathbf{B}^\ast$ satisfies the following matrix equation:
\begin{align}\label{B_solution}\mathbf{A} \mathbf{A}^\top \mathbf{B} \mathbf{\Sigma}_{\mathbf{yy}} + \lambda \nabla \rho(\mathbf{B}) = \mathbf{A} \mathbf{\Sigma}_{\mathbf{xy}}.
\end{align}
with $\mathbf{\Sigma}_{\mathbf{xy}}
    = \mathbb{E}\mathbf{x} \mathbf{y}^\top$.

Having established existence and uniqueness of the learned filter, we now analyze the stability of the solution with respect to perturbations in the data distribution. Let $\mathbf{B}^*_\pi$ and $\mathbf{B}^*_{\pi'}$ be the minimizers of (\ref{Learned Inverse_supervised}) under distributions $(\mathbf{x},\mathbf{y})\sim{\pi}$ and $(\mathbf{x},\mathbf{y})\sim{\pi'}$, respectively. Given Assumptions (A1)-(A3), there exists a constant $C > 0$ such that:$$\|\mathbf{B}^*_{\pi} - \mathbf{B}^*_{\pi'}\|_F \leq \frac{C}{\lambda \mu} \, W_2(\pi, \pi'),$$
where $W_2$ denotes the 2-Wasserstein distance. A proof is given in Appendix \ref{Stability_proof}. 





\subsubsection{Interpretation as a regularized pseudo-inverse}\label{Interpretation}
A more detailed interpretation of the learned filter can be obtained by considering a concrete example. In addition to assumptions (A1)-(A3), we let
\begin{itemize}
    \item Additive Gaussian noise: $\mathbf{y} = \mathbf{A}\mathbf{x} + \mathbf{n}$ with $\mathbf{n}\sim \mathcal{N}(\mathbf{0},\sigma^2 \mathbf{I})$,
    \item Tikhonov regularization: $\rho(\mathbf{B}) = \|\mathbf{B}\|_F^2$.
\end{itemize}
We then have
\begin{equation}\label{statistics_def}
\mathbf{\Sigma}_{\mathbf{yy}} = \mathbf{A} \mathbf{\Sigma}_{\mathbf{xx}} \mathbf{A}^\top+\sigma^2\mathbf{I}, \qquad \mathbf{\Sigma}_{\mathbf{xy}} = \mathbf{\Sigma}_{\mathbf{xx}} \mathbf{A}^\top.
\end{equation}

Combining (\ref{B_solution}) and (\ref{statistics_def}) we then find that $\mathbf{B}^\ast$ satisfies the following matrix equation
\begin{equation}
\mathbf{AA}^\top\mathbf{B}\mathbf{A}\mathbf{\Sigma}_{\mathbf{xx}}\mathbf{A}^\top+\sigma^2\mathbf{AA}^\top\mathbf{B}+\lambda\mathbf{B}=\mathbf{A}\mathbf{\Sigma}_{\mathbf{xx}}\mathbf{A}^\top.
    \label{eq:kernel_tr_form}
\end{equation}
Introducing the auxiliary matrices
\begin{equation}
    \mathbf{M} = \mathbf{A}\mathbf{A}^\top,
    \;
    \mathbf{N} =\mathbf{A}\mathbf{\Sigma}_{\mathbf{xx}}\mathbf{A}^\top,
\end{equation}
we can re-write this as
\begin{equation}
    (\mathbf{N}\otimes \mathbf{M} + \sigma^2 \mathbf{I}\otimes\mathbf{M}+\lambda\mathbf{I} )\,\mathrm{vec}(\mathbf{B}) = \mathrm{vec}(\mathbf{N}),
    \label{eq:sylvester}
\end{equation}
We can immediately see that when $\sigma = 0$, $\lambda =0$ we learn the regular pseudo-inverse: $\mathbf{B}^* = (\mathbf{A}\mathbf{A}^\top)^{-1}$.

To further interpret the learned kernel and its relation to classical analytical filters, we analyze the structure of~\eqref{eq:sylvester} in the spectral domain.  
Let the singular value decomposition (SVD) of the forward operator be $\mathbf{A} = \mathbf{U}{\mathbf{\Sigma}}\mathbf{V}^\top$.  
Substituting this decomposition into~\eqref{eq:sylvester} and premultiplying and postmultiplying by~$\mathbf{U}^\top$ and~$\mathbf{U}$, respectively, yields
\begin{equation}
    \boldsymbol{\Lambda}\,\mathbf{H}\,\widetilde{\mathbf{N}} + (\sigma^2\mathbf{\Lambda}+\lambda\mathbf{I}) \mathbf{H} = \widetilde{\mathbf{N}},
    \label{eq:Btilde}
\end{equation}
where
\[
\boldsymbol{\Lambda} = \operatorname{diag}(\sigma_1^2,\dots,{\sigma}_M^2), \quad
\mathbf{H} = \mathbf{U}^\top \mathbf{B}\mathbf{U}, \quad
\widetilde{\mathbf{N}} = \mathbf{U}^\top (\mathbf{A}\mathbf{\Sigma_{xx}}\mathbf{A}^\top) \mathbf{U}.
\]
If in addition we assume $\mathbf{\Sigma}_{\mathbf{xx}}=\mathbf{I}$, $\widetilde{\mathbf{N}}$ is diagonal and the optimal $\mathbf{H}$ is diagonal as well. In this case, the spectral components decouple, and~\eqref{eq:Btilde} reduces entry-wise to
\begin{equation}
    {\sigma}_i^4 H_{ii}\, + (\sigma^2{\sigma}_i^2+\lambda) H_{ii} = \sigma_i^2,
\end{equation}
leading to the closed-form solution
\begin{equation}
    H_{ii}
    = \frac{\sigma_i^2}{(\sigma_i^2 + \sigma^2){\sigma}_i^2+\lambda}.
    \label{eq:b_diag}
\end{equation}
Each $H_{ii}$ thus acts as a spectral gain that modulates the contribution of the $i$-th mode according to its singular value, and noise variance. It provides an interpretable, data-adaptive filtering mechanism.  
Modes with small~${\sigma}_i$ are damped by the regularization term~$\lambda$, effectively suppressing noise amplification. This interpretation highlights that the learned kernel $\mathbf{B}$ behaves as a \emph{regularized pseudo inverse} in the spectral domain, adapting its frequency response based on the data and noise statistics. This expression closely resembles classical frequency-domain filters, such as the ramp filter in filtered back projection (FBP) but with an additional data-adaptive weighting that depends on the signal and noise variances.


\subsubsection{Reconstruction error}
Given a solution $\mathbf{B}^*$ to \eqref{Learned Inverse_supervised}, and under the same assumptions as in section 3.1.3, we can decompose the resulting reconstruction error as (see Appendix \ref{Reconstruction_error})

\begin{align}
\label{rec_err_bound}
    \mathbb{E}
    \|\mathbf{A}^\top \mathbf{B}^* \mathbf{y} - \mathbf{x}\|_2^2
    &\le
    \mathbb{E}
    \|\mathbf{A}^\top \mathbf{B}^* \mathbf{n}\|_2^2
    \\\nonumber
    &\quad +
    \mathbb{E}
    \|(\mathbf{A}^\top \mathbf{B}^* - \mathbf{A}^\dagger) \mathbf{Ax}\|_2^2
    \\\nonumber
    &\quad +
    \mathbb{E}
    \|(\mathbf{A}^\dagger \mathbf{A} - \mathbf{I})\mathbf{x}\|_2^2.
\end{align}
In the upper bound, the first term is interpreted as the variance (the influence of noise on the
reconstruction) while the second term is interpreted as the bias (the influence of regularization on the learned inverse). The last term indicates a null-space error that we cannot get rid of, as the learned inverse strictly operates in the column space of $\mathbf{A}$. In the following, we investigate the first two error terms in the spectral domain. 

We begin from the first term of the upper bound in (\ref{rec_err_bound}). Using decompositions introduced for $\mathbf{A}$ and $\mathbf{B}$ in section (\ref{Interpretation}), we can write $\mathbf{A^{\top}B}$ as,
\begin{equation}
    \mathbf{A^\top B}= \mathbf{V\Sigma HU^{\top}}=\mathbf{V}\,\mathrm{diag}\left( \frac{\sigma_i^3}{(\sigma_i^2+ \sigma^2)\sigma_i^2 + \lambda} \right)\mathbf{U}^{\top},
\end{equation}
where, we define the spectral variance factor as,
\begin{align}
    \zeta_i=\frac{\sigma_i^3}{(\sigma_i^2+ \sigma^2)\sigma_i^2 + \lambda},
\end{align}
then, it can be written as,
\begin{equation}
       \mathbb{E}\|\mathbf{A}^\top \mathbf{B} \mathbf{n}\|_2^2=\sigma^2\cdot\sum_{i=1}^M  \zeta_i^2
\end{equation}
Above derivation shows that the variance has a direct relation with the noise variance. Moreover, it goes to zero as the regularization parameter $\lambda$ becomes very large.

The bias term of the upper bound in (\ref{rec_err_bound}) can be investigated in the spectral domain as well. To do so, we start with:
\begin{align}
(\mathbf{A}^{\top}\mathbf{B} - \mathbf{A}^{\dagger})\mathbf{A}=\Big(\mathbf{V\Sigma H}\mathbf{U}^{\top} - \mathbf{V}\,\mathbf{\Sigma^{\dagger}}
\mathbf{U}^{\top}\Big)\mathbf{U\Sigma V^{\top}} = \mathbf{V}\,\mathrm{diag}\Big(\frac{\sigma_{i}^{4}}{(\sigma_{i}^{2} + \sigma^{2})\sigma_{i}^{2} + \lambda}-1\Big)\mathbf{V}^\top
\end{align} 
 Then, we define the spectral bias factor as,
\begin{align}
    \beta_i = 1-\frac{\sigma_{i}^{4}}{(\sigma_{i}^{2} + \sigma^{2})\sigma_{i}^{2} + \lambda}=\frac{\sigma^{2}\sigma_{i}^{2} + \lambda}
{\sigma_{i}^{4} + \sigma^{2}\sigma_{i}^{2} + \lambda},
\end{align}
 since we assume $\mathbf{x}\sim\mathcal{N}(0,\mathbf{I})$, then it reads:
\begin{align}
    \mathbb{E}\!\left\|
    (\mathbf{A}^{\top}\mathbf{B} - \mathbf{A}^{\dagger})\mathbf{A}\mathbf{x}
    \right\|_{2}^{2}
    =
    \sum_{i=1}^{M}
    \beta_i^2.
\end{align}
It indicates that the bias term vanishes for very large singular values. Furthermore, for very small $\sigma_i$'s, it goes to $1$ and introduces bias error.


Finally, we have an error-component that we cannot get rid of as the learned inverse is fundamentally restricted to producing solutions in the range of $\mathbf{A}^\top$ and thus cannot fill in components in the null-space.


\subsection{Parameterizing the filter for tomographic reconstruction}
For computational efficiency we want to learn $\mathbf{B}$ as a composition of filtering and weighting steps, parametrized by a set of parameters $\mathbf{p}$. The form of $\mathbf{B}(\mathbf{p})$ can be informed by the acquisition geometry, as we will see below. The choice for $\rho$ will be discussed in more detail as well.

The well-posedness results in the previous section are proved for the unconstrained problem over the full matrix space,
i.e., $\mathbf{B}\in\mathbb{R}^{M\times M}$ with $\lambda>0$.
When introducing a parametrization $\mathbf{B}=\mathbf{B}(\mathbf{p})$, we restrict the admissible operators to
\[
\mathcal{S}=\{\mathbf{B}(\mathbf{p}) : \mathbf{p}\in\mathcal{P}\},
\]
and consequently solve a constrained problem over $\mathcal{S}$ rather than over all of $\mathbb{R}^{M\times M}$.

\smallskip
\noindent\emph{Existence.}
Existence of a minimizer typically carries over under mild assumptions, e.g., if
$\mathbf{p}\mapsto \mathbf{B}(\mathbf{p})$ is continuous and $\mathcal{P}$ is compact (or the objective is coercive in
$\mathbf{p}$), then the continuous objective function attains its minimum on a compact set \cite{Bolzano}.

\smallskip
\noindent\emph{Uniqueness and stability.}
Uniqueness and stability do not automatically carry over in parameter space.
Even if the minimizer is unique in operator space (i.e., the best $\mathbf{B}$ within $\mathcal{S}$ is unique),
the corresponding parameter $\mathbf{p}$ may be non-unique unless the parametrization is identifiable \cite{Petrica2024ANO}.
Moreover, if $\mathbf{B}(\mathbf{p})$ depends nonlinearly on $\mathbf{p}$ (e.g., through bilinear compositions such as
\emph{weights} $\times$ \emph{filtering}), then the objective is generally nonconvex in $\mathbf{p}$, so global uniqueness
is typically lost and only local guarantees can be expected.

\smallskip
Independently of well-posedness, restricting to $\mathcal{S}$ introduces an approximation aspect: the parametrized solution
is the best admissible operator in $\mathcal{S}$, which coincides with the unconstrained optimizer $\mathbf{B}^\ast$
only if $\mathbf{B}^\ast\in\mathcal{S}$ (or can be well-approximated by elements of $\mathcal{S}$).

\smallskip
\noindent In summary, the well-posedness results of Section~\ref{wellposedness} remain a useful reference at the operator level
$\mathbf{B}$. In the parametrized setting, existence is typically preserved under mild conditions, whereas uniqueness and
stability in the parameters depend on structural properties of the map $\mathbf{p}\mapsto \mathbf{B}(\mathbf{p})$,
for instance whether it preserves convexity (e.g., affine dependence on $\mathbf{p}$) and whether it is identifiable.

\smallskip
\noindent Below, we specify geometry-informed parameterizations of $\mathbf{B}(\mathbf{p})$ for the acquisition settings considered in this work.
\subsubsection{2D Parallel beam} 
For the parallel-beam geometry with $m$ angles and a detector with $d$ pixels (so $M = m\cdot d$), $\mathbf{B}$ is expressed as the 1D convolution along the detector-pixel dimension with a specified filter. Specifically, we can present it as
\begin{equation}\label{FBP_Parallel}
    \mathbf{B}(\mathbf{p})
    = \mathbf{I}_m \otimes \big(\mathbf{F}_d^{-1}\operatorname{diag}(\mathbf{p})\mathbf{F}_d \big),
\end{equation}
where $\mathbf{I}_m$ represents the $m\times m$ identity matrix, $\mathbf{F}_d$ the $d-$dimensional DFT and $\mathbf{p}\in\mathbb{R}^{d}$ represent the filter coefficients.
As analytically derived in the literature~\cite{Hensen}, in the ideal continuous and noise-free scenario, $\mathbf{p}$ corresponds to a high-pass ramp filter.
\subsubsection{2D Circular fan-beam}
In a circular fan-beam geometry with $m$ angles and a detector with $d$ pixels (so $M = m\cdot d$), an additional weighting factor is needed to compensate for the non-uniform ray coverage \cite{Kak}. We associate the learnable parameter $\mathbf{p}$ with a filter $\mathbf{p}_\text{filter} \in \mathbb{R}^d$ and weight $\mathbf{p}_\text{weight}\in\mathbb{R}^{d}$. We then have
\begin{equation}\label{FBP_Fan}
\mathbf{B}(\mathbf{p}) = \mathbf{I}_m \otimes \big(\mathbf{F}_d^{-1} \operatorname{diag}(\mathbf{p}_\text{filter})\mathbf{F}_d\operatorname{diag}(\mathbf{p}_\text{weight})\big)
\end{equation}
\subsubsection{2D Elliptical fan-beam} This acquisition scheme, where the source and detector rotate around the object along an elliptical trajectory, is commonly used to accommodate objects with large aspect ratios. As the geometry changes per projection angle, the kernel $\mathbf{B}$ requires a parametrization for the filter and weighting vector that is slightly different from the circular fan beam case. Namely, we consider angle-dependent filters and weights and represent $\mathbf{B}$ as below,
\begin{equation}
\label{F_recon:elips}
\mathbf{B}(\mathbf{p}) = \big(\mathbf{I}_m\otimes \mathbf{F}_d^{-1}\big)\operatorname{diag}(\mathbf{p}_\text{filter})\big(\mathbf{I}_m\otimes \mathbf{F}_d\big)\operatorname{diag}(\mathbf{p}_\text{weight}),
\end{equation}
where $\mathbf{p}_\text{filter}\in\mathbb{R}^{M}$, $\mathbf{p}_\text{weight}\in\mathbb{R}^{M}$ represent an angle-dependent filter and weight.
\subsubsection{3D Circular laminographic trajectory:}
\begin{figure}
    \centering
\includegraphics[width=0.5\linewidth]{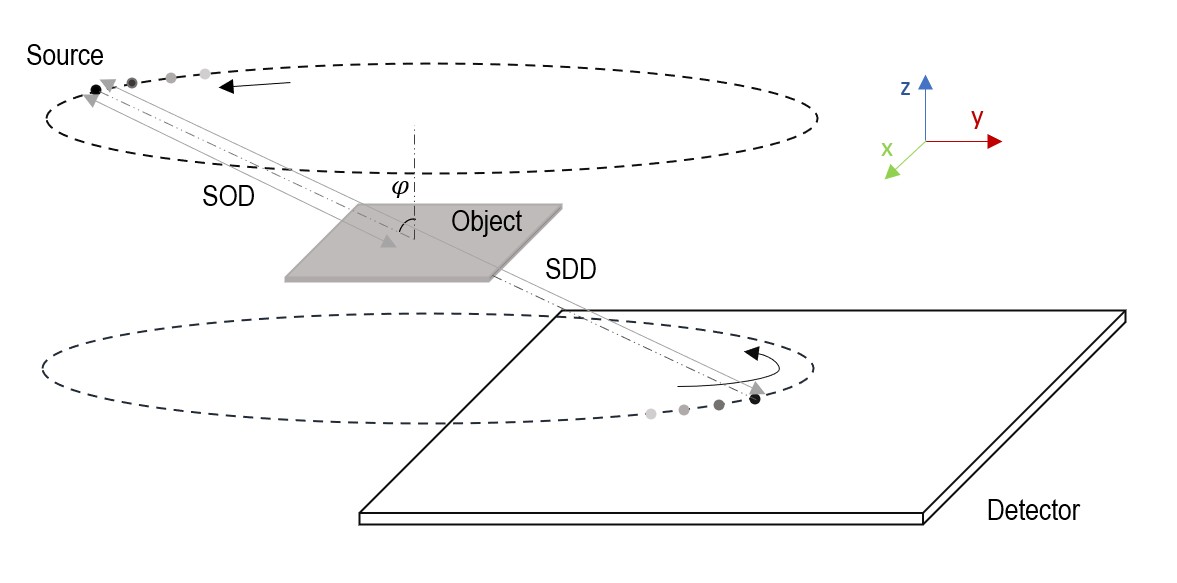}
    \caption{Sketch of acquisition geometry in a circular laminography setup.}
    \label{fig:Lamino-setup}
\end{figure}

Laminography is an imaging technique tailored for flat, layered objects. Several acquisition trajectories are used in practice depending on sample size and hardware constraints \cite{Brien}. Here we focus on the \emph{circular} laminography trajectory, which is well suited for large, flat samples: a point source and a flat detector rotate on opposite sides of the sample with a fixed laminography tilt angle~$\phi$ (the central ray is tilted by $\phi$ with respect to the sample plane); see Fig.~\ref{fig:Lamino-setup}. Throughout the rotation, the central ray is aligned to pass through a fixed point in the sample (typically the center of the middle layer), ensuring a consistent geometry.

In this setup, the detector data per projection angle are two-dimensional, hence we parametrize the kernel $\mathbf{B}$ using angle-dependent 2D filters and 2D weights. Let $\mathbf{F}_{d_1}$ and $\mathbf{F}_{d_2}$ denote the 1D Fourier transforms acting along the two detector dimensions. Exploiting the separability of the 2D Fourier transform,
\[
\mathbf{F}_{d_1\cdot d_2} = \mathbf{F}_{d_1}\otimes \mathbf{F}_{d_2},
\]
the 2D Fourier transform acting independently on each view can be written as a Kronecker product of 1D Fourier operators.

Therefore, we parametrize $\mathbf{B}$ as
\begin{equation}
\label{F_recon:circular}
\mathbf{B}(\mathbf{p})
=
\big(\mathbf{I}_m \otimes \mathbf{F}_{d_1}^{-1} \otimes \mathbf{F}_{d_2}^{-1}\big)
\operatorname{diag}(\mathbf{p}_{\text{filter}})
\big(\mathbf{I}_m \otimes \mathbf{F}_{d_1} \otimes \mathbf{F}_{d_2}\big)
\operatorname{diag}(\mathbf{p}_{\text{weight}}),
\end{equation}
where $\mathbf{p}_{\text{filter}}\in\mathbb{R}^{M}$ and $\mathbf{p}_{\text{weight}}\in\mathbb{R}^{M}$ are vectors of concatenated 2D filter and 2D weight coefficients, with $M = m\cdot d_1\cdot d_2$.

\subsection{Implementation}
For given training set $\{(\mathbf{x}_i,\mathbf{y}_i)\}_{i=1}^K$ we now need to solve 
\begin{align}\label{supervised_impl}
\min_{\mathbf{p}} K^{-1}\sum_{i=1}^K
    \|\mathbf{A}^\top \mathbf{B}(\mathbf{p}) \mathbf{y}_i - \mathbf{x}_i\|_2^2 + \lambda \rho(\mathbf{p}),    
\end{align}

where $\mathbf{p}=(\mathbf{p}_\text{filter}, \mathbf{p}_\text{weight})$ contains the parameters for the filter and weight. In contrast to~\eqref{Learned Inverse_supervised}, the parametrized formulation is, in general, no longer strictly convex. While the objective is strictly convex in the operator $\mathbf{B}$, substituting $\mathbf{B}=\mathbf{B}(\mathbf{p})$ changes the geometry of the problem: if $\mathbf{B}(\mathbf{p})$ depends non-affinely on $\mathbf{p}$, i.e., $\mathbf{B}$ is formed by composition or products of parameters, then the resulting objective becomes nonconvex in $\mathbf{p}$. Consequently, global uniqueness in parameter space is not guaranteed and optimization typically converges to a first-order critical point of the objective. The resulting learned operator therefore corresponds to a (possibly local) solution within the admissible family $\mathcal{S}=\{\mathbf{B}(\mathbf{p})\}$.

\subsubsection{Regularization}
We generally want to enforce smoothness on both the filter and spatial weights, which we achieve by setting
\[
\rho(\mathbf{p}) = \|\mathbf{L}_\text{freq}\mathbf{p}_\text{filter}\|_2^2 + \|\mathbf{L}_\text{space}\mathbf{p}_\text{weight}\|_2^2,
\]
with $\mathbf{L}_\text{freq}, \mathbf{L}_\text{space}$ representing a finite-difference approximation of the derivative or gradient along the frequency and spatial axes. In addition, we enforce symmetry of the filter, which guarantees that its Fourier transform is real-valued.

\subsubsection{Implementation Details}
The forward operator $\mathbf{A}$ is implemented using \text{Tomosipo} \cite{hendriksen-2021-tomos}, providing a \text{PyTorch}-compatible interface to the \text{ASTRA} toolbox \cite{Astra}. For the training process, we utilize \text{PyTorch}'s automatic differentiation to perform backpropagation through the operator $\mathbf{B}(\mathbf{p})$. The entire pipeline — from the frequency-domain parameters to the final reconstruction error — is differentiable, allowing us to utilize a gradient-based optimizer (Adam) to update both the filter coefficients and spatial weights simultaneously \cite{Kingma2014AdamAM}.
\section{Case studies}

In this section, we investigate the data-driven filter design method under different practical scenarios. We begin by learning parameters for noisy measurements in both 2D parallel-beam and 2D fan-beam geometries. Later, we learn filters and weights for a 2D elliptical trajectory and a 3D circular laminography setup and show numerical results.

\subsection{2D Parallel Beam Geometry}


We focus on learning a filter for noisy measurements in parallel-beam geometry. Specifically, we investigate how the filter in (\ref{FBP_Parallel}) can be adapted to noisy scenarios by making use of training data. To this end, we generate synthetic 2D phantoms of size $400\times400$ consisting of circles with random sizes and positions with binary values. The simulated acquisition consists of $360$ uniformly spaced projections with a resolution of $512$ pixels. The pixel size is set to $0.002$ for both the phantom and the detector. To ensure stability during training, we include a small smoothness regularization with the coefficient of $0.001$ to suppress abrupt variations in the filter. 

\begin{figure}
    \centering
    \includegraphics[width=0.5\linewidth]{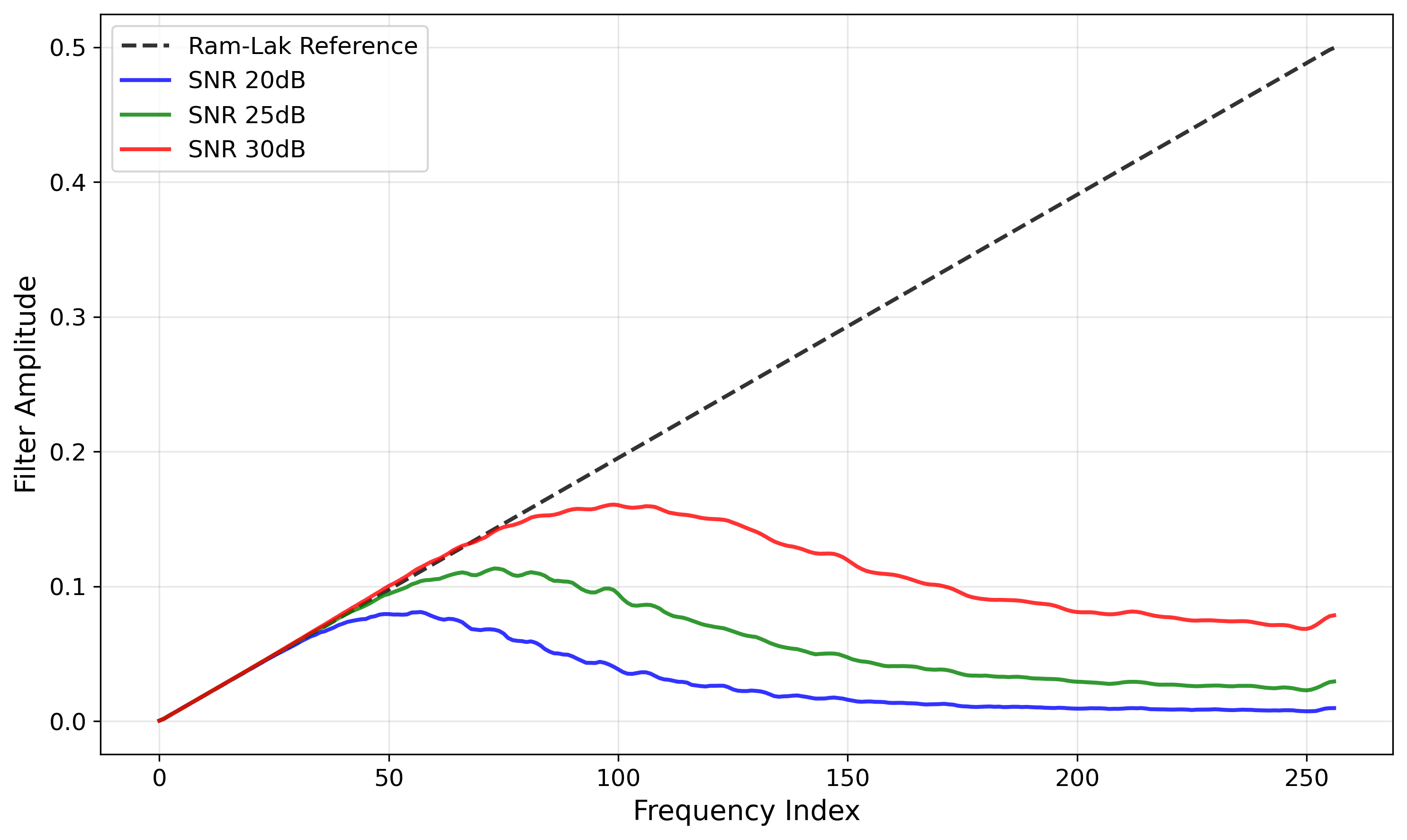}
    \caption{Learned filters for different noisy datasets in a parallel beam geometry versus Ram-Lak filter shown in the frequency domain.}
    \label{filter_parallel}
\end{figure}

Figure \ref{filter_parallel} shows the filters learned from datasets corrupted with different levels of additive Gaussian noise, along with the standard Ram-Lak filter, which is typically used in FBP. While the Ram-Lak filter exhibits a linear increase in amplitude with frequency,  the learned filters follow a similar trend at low frequencies but begin to deviate and attenuate at higher frequencies across all noise scenarios, with the inflection point shifting toward lower frequencies for higher noise levels. This behavior suggests that the learned filters adapt to the noise characteristics of the training data, effectively acting as modified high-pass filters that suppress noise-induced artifacts while the unmodified Ram-Lak filter might otherwise amplify. Notably, the same behavior with respect to noise is commonly observed or enforced in other adaptive and fixed filter types \cite{Dan, Shepp_Logan_fil}.

We evaluate the performance of the learned filters on separate validation datasets with the same noise levels. Table~\ref{tab:metrics_parallel} presents the MSE and SSIM results across various noise conditions. From these results, we observe that the standard FBP method becomes less effective at higher noise levels, while at lower noise levels its performance becomes comparable with learned method. Moreover, the MSE results in Table~\ref{tab:metrics_parallel} are the best when the validation noise level matches the training noise level, although this trend does not hold consistently for SSIM. In particular, SSIM may favor reconstructions that are smoother and contain fewer high-frequency fluctuations, whereas MSE penalizes bias and oversmoothing more directly. Consequently, a model trained at a different noise level can achieve higher SSIM on a given test set if it suppresses high-frequency artifacts even when its MSE is not optimal.


\begin{table}[ht]
\centering
\caption{Learned reconstruction performance for different training/validation SNR levels in parallel-beam geometry (T/V: Train SNR / Val SNR). Entries are reported as \textbf{MSE} (mean $\pm$ std) $\times 10^{-3}$ \,/\, \textbf{SSIM} (mean $\pm$ std).}
\label{tab:metrics_parallel}
\setlength{\tabcolsep}{6pt}
\renewcommand{\arraystretch}{1.15}
\begin{tabular}{lccc}
\toprule
\textbf{T/V} & \textbf{20} & \textbf{25} & \textbf{30} \\
\midrule
\textbf{20} &
$5.3 \pm 0.6$ / $0.305 \pm 0.034$ &
$4.3 \pm 0.4$ / $0.499 \pm 0.035$ &
$4.0 \pm 0.3$ / $0.646 \pm 0.031$ \\
\textbf{25} &
$6.2 \pm 1.0$ / $0.197 \pm 0.023$ &
$3.9 \pm 0.4$ / $0.369 \pm 0.037$ &
$3.2 \pm 0.3$ / $0.564 \pm 0.039$ \\
\textbf{30} &
$10.1 \pm 2.1$ / $0.132 \pm 0.0098$ &
$4.6 \pm 0.7$ / $0.247 \pm 0.026$ &
$2.9 \pm 0.3$ / $0.444 \pm 0.043$ \\
\midrule
\textbf{FBP} &
$39 \pm 9$ / $0.080 \pm 0.001$ &
$12 \pm 2$ / $0.120 \pm 0.005$ &
$4 \pm 1$ / $0.220 \pm 0.022$ \\
\bottomrule
\end{tabular}
\end{table}

\clearpage
\subsection{2D Fan Beam Geometry}
This section extends the learning of parameters (filter and weights) for the FBP approach to fan-beam geometry under noisy conditions. To learn the required FBP parameters using clean training data, it is sufficient to incorporate equation~(\ref{FBP_Fan}) into the optimization problem of~(\ref{supervised_impl}). The training and validation datasets consist of synthetic 2D phantoms of size $400\times400$ containing circles with random sizes and positions with pixel size of $0.005$. The simulated acquisition includes $360$ projections with a resolution of $1024$ pixels and a pixel size of $0.01$, while the source-to-object distance ($SOD$) and source-to-detector distance ($SDD$) are set to $2.5$ and $5$, respectively. To ensure stable training, a small smoothness regularization term with a coefficient of $0.0001$ is added to independently suppress abrupt variations in both the filter and the weight vector.

\begin{figure}
    \centering
    \begin{subfigure}{0.45\textwidth}
         \centering
         \includegraphics[width=\linewidth]{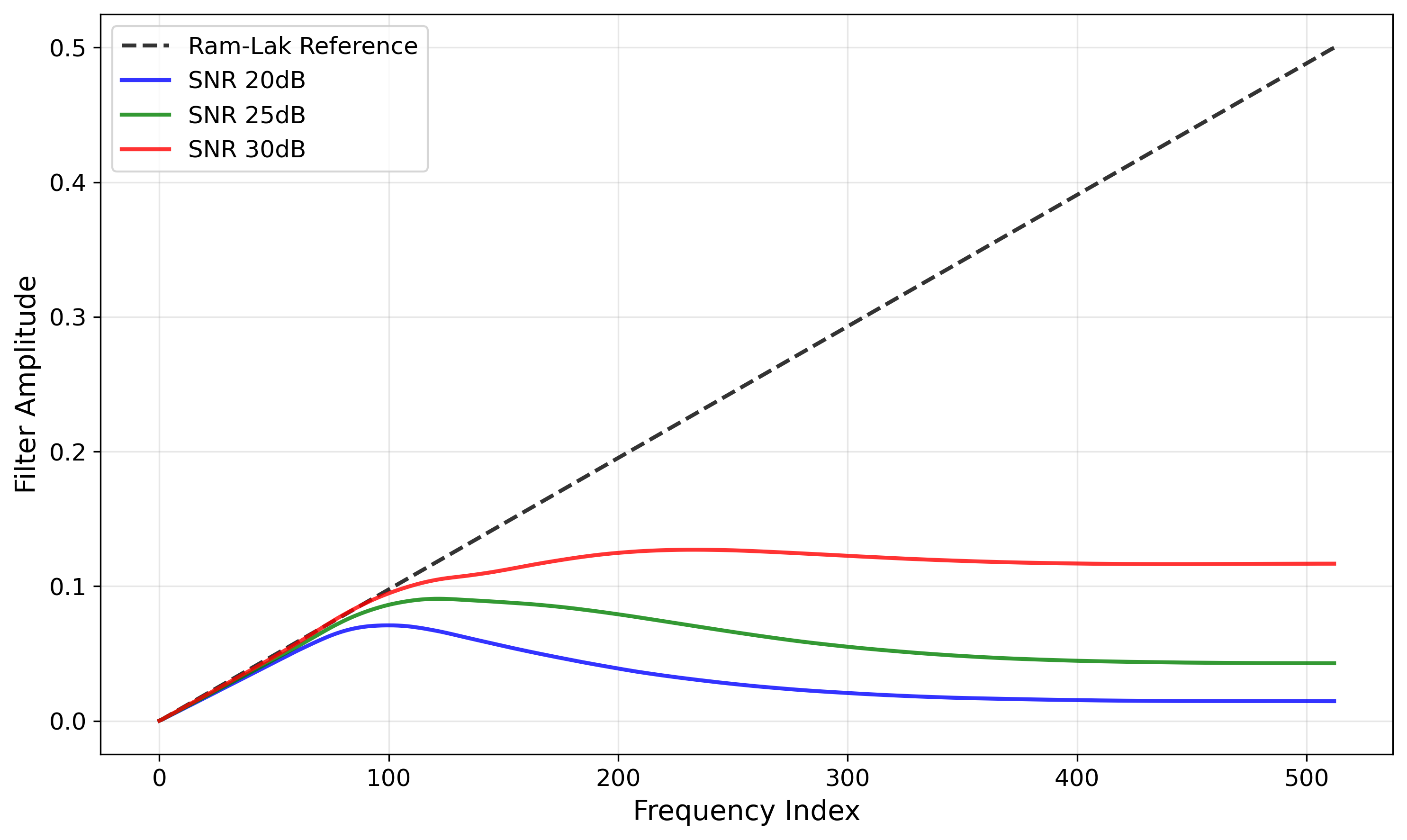}
         \caption{}
    \end{subfigure}
    \begin{subfigure}{0.45\textwidth}
        \centering
        \includegraphics[width=\linewidth]{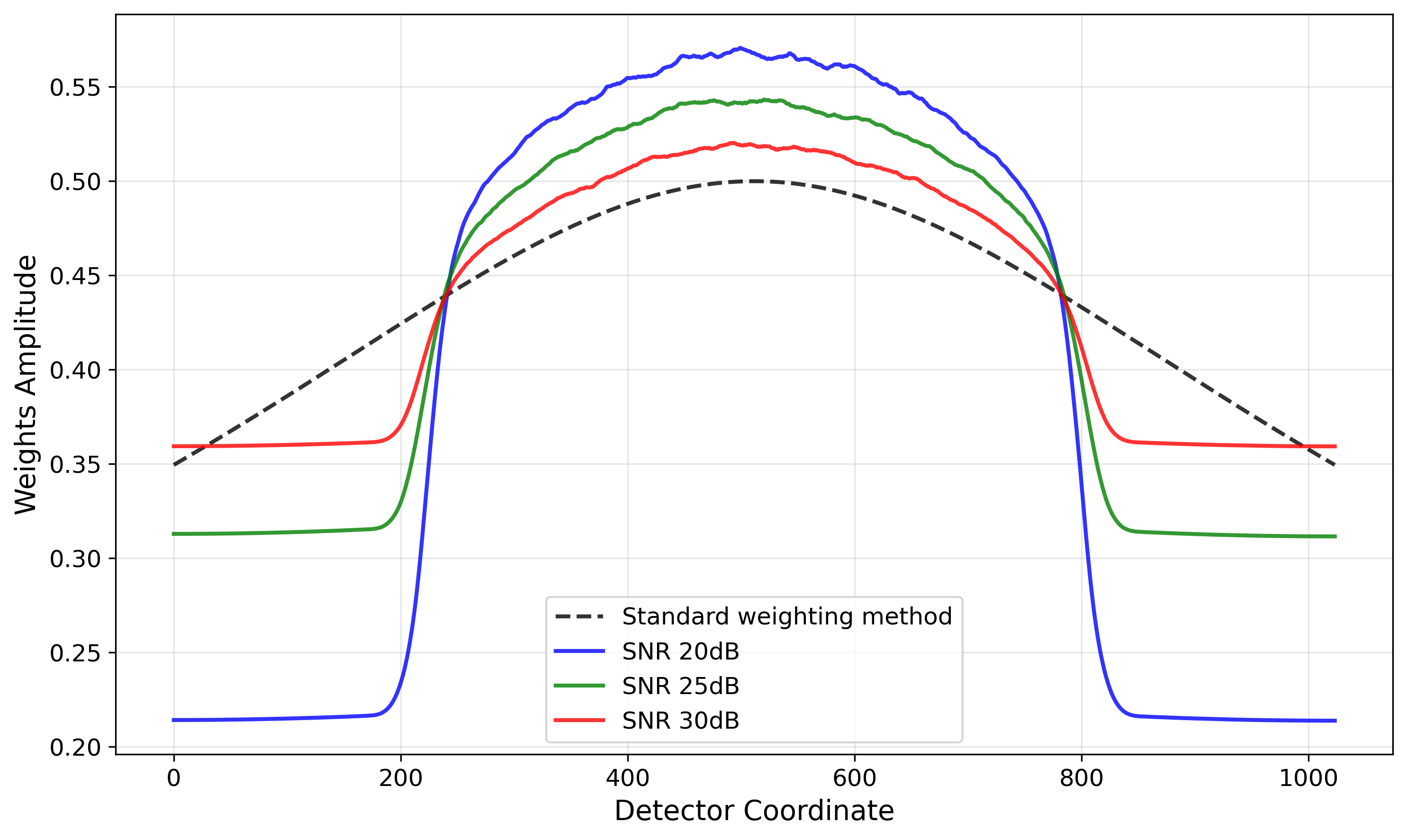}
        \caption{}
    \end{subfigure}
    \caption{Learned fan-beam (A) filter and (B) weights for different noise levels (SNR 20/25/30 dB), compared to the Ram-Lak filter and the standard fan-beam weighting used in conventional FBP. The weights are normalized to the standard profile to remove the filter--weight scaling ambiguity; the filter is rescaled accordingly.}
    \label{filter:weights:fanbeam}
\end{figure}

Figure~\ref{filter:weights:fanbeam} illustrates the learned filter and fan-beam weighting for different noise levels, together with the standard choices used in conventional FBP. As in the parallel-beam case, the learned filters deviate from the Ram-Lak at higher frequencies, with stronger attenuation at lower SNR, reflecting noise-adaptive spectral shaping. Moreover, the learned weights mimic the standard weighting method with noise-dependent deviations that are largest at low SNR. Since the fan-beam parametrization is factored into a filter and weights, the representation is subject to a global scale ambiguity (a rescaling of the filter can be compensated by an inverse rescaling of the weights). Following common practice in bilinear/factored models, we fix this gauge by normalizing the learned weights to the standard weighting profile and rescaling the learned filter accordingly \cite{Jin_2018_ECCV}.

Quantitative results for the fan-beam case, presented in Table~\ref{tab:metrics_fan}, further confirm the advantage of the learned filters over the FBP method. While FBP performs poorly across all noise levels, the learned models consistently achieve lower MSE and higher SSIM values. As in the parallel-beam experiments, the best MSE results occur when the training and validation noise levels are matched. However, this pattern is not always reflected in SSIM, similarly to the parallel beam case.

\begin{table}[ht]
\centering
\caption{Learned reconstruction performance for different training/validation SNR levels in fan-beam geometry (T/V: Train SNR / Val SNR). Entries are reported as \textbf{MSE} (mean $\pm$ std) $\times 10^{-3}$ \,/\, \textbf{SSIM} (mean $\pm$ std).}
\label{tab:metrics_fan}
\setlength{\tabcolsep}{6pt}
\renewcommand{\arraystretch}{1.15}
\begin{tabular}{lccc}
\toprule
\textbf{T/V} & \textbf{20} & \textbf{25} & \textbf{30} \\
\midrule
\textbf{20} &
$5.0 \pm 0.5$ / $0.319 \pm 0.028$ &
$4.1 \pm 0.4$ / $0.504 \pm 0.041$ &
$3.8 \pm 0.3$ / $0.629 \pm 0.027$ \\
\textbf{25} &
$6.1 \pm 0.8$ / $0.217 \pm 0.019$ &
$4.0 \pm 0.5$ / $0.381 \pm 0.040$ &
$3.4 \pm 0.3$ / $0.525 \pm 0.030$ \\
\textbf{30} &
$9.9 \pm 1.7$ / $0.147 \pm 0.009$ &
$4.9 \pm 0.8$ / $0.262 \pm 0.031$ &
$3.4 \pm 0.4$ / $0.401 \pm 0.031$ \\
\midrule
\textbf{FBP} &
$38.5 \pm 7.7$ / $0.090 \pm 0.001$ &
$13.8 \pm 3.2$ / $0.133 \pm 0.008$ &
$6.3 \pm 1.1$ / $0.205 \pm 0.016$ \\
\bottomrule
\end{tabular}
\end{table}

\clearpage
\subsection{2D Elliptical Trajectory with Fan Beam}
In this case, we learn filter and weights for filtered back-projection in a 2D elliptical trajectory. It requires incorporating~(\ref{F_recon:elips}) into~(\ref{supervised_impl}) and minimize over the filter and weights. The training dataset includes synthetic 2D phantoms of size $400\times400$ with random circles in terms of size and position. The measurements are acquired from 360 projections over a full rotation, each consisting of 512 pixels with a pixel size of 0.015. The source and detector rotate around the object along an elliptical trajectory, where the source-to-object distance (SOD) and source-to-detector distance (SDD) are 4.5 and 9, respectively, along the semi-minor axis, and 12 and 24 along the semi-major axis. 

\begin{figure}[htbp]
    \centering
    


    \begin{subfigure}{0.32\textwidth}
        \centering
        \includegraphics[width=\linewidth]{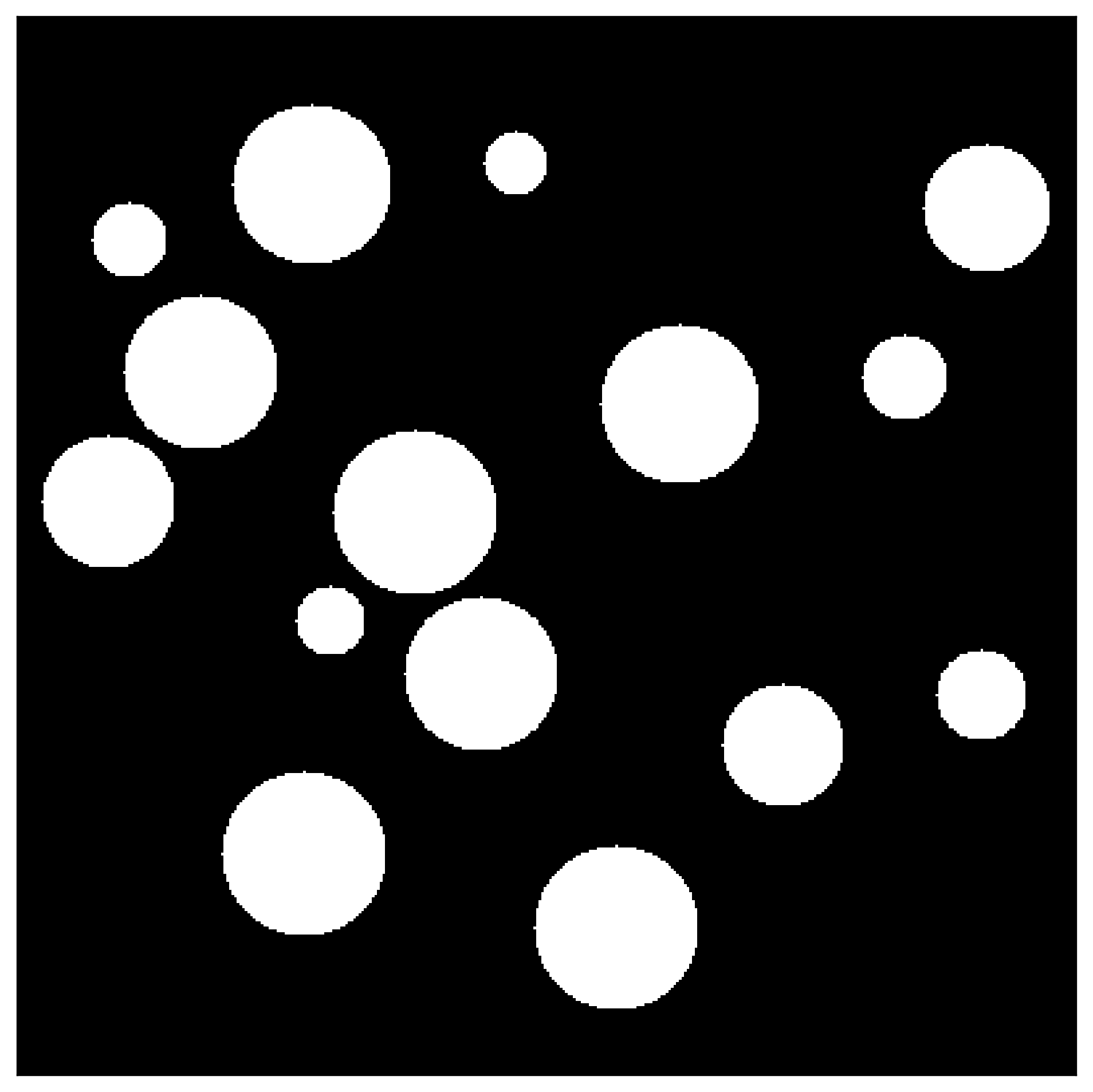}
        \caption*{Ground Truth}
    \end{subfigure}\hfill
    \begin{subfigure}{0.32\textwidth}
        \centering
        \includegraphics[width=\linewidth]{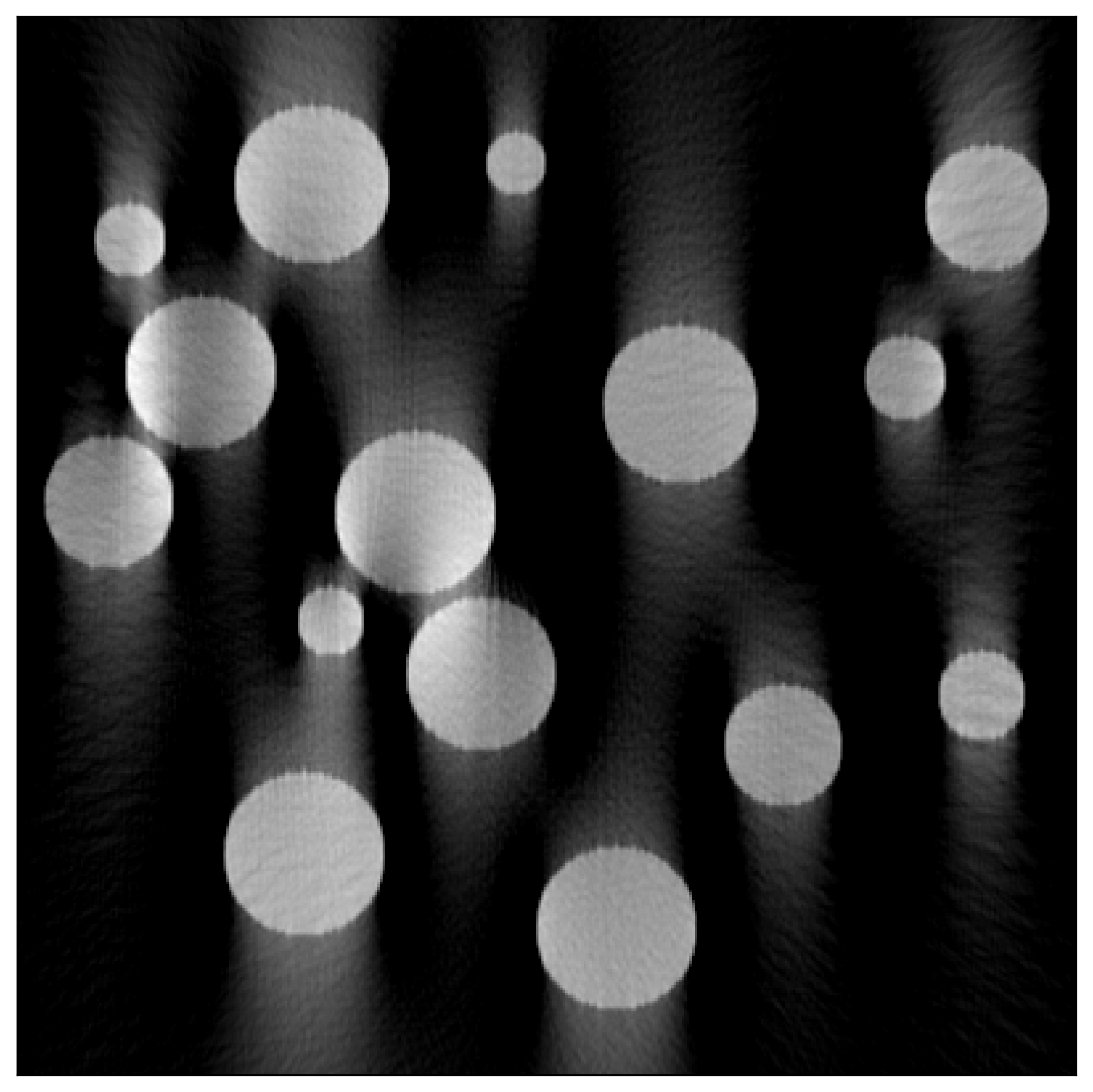}
        \caption*{Conventional FBP}
    \end{subfigure}\hfill
    \begin{subfigure}{0.32\textwidth}
        \centering
        \includegraphics[width=\linewidth]{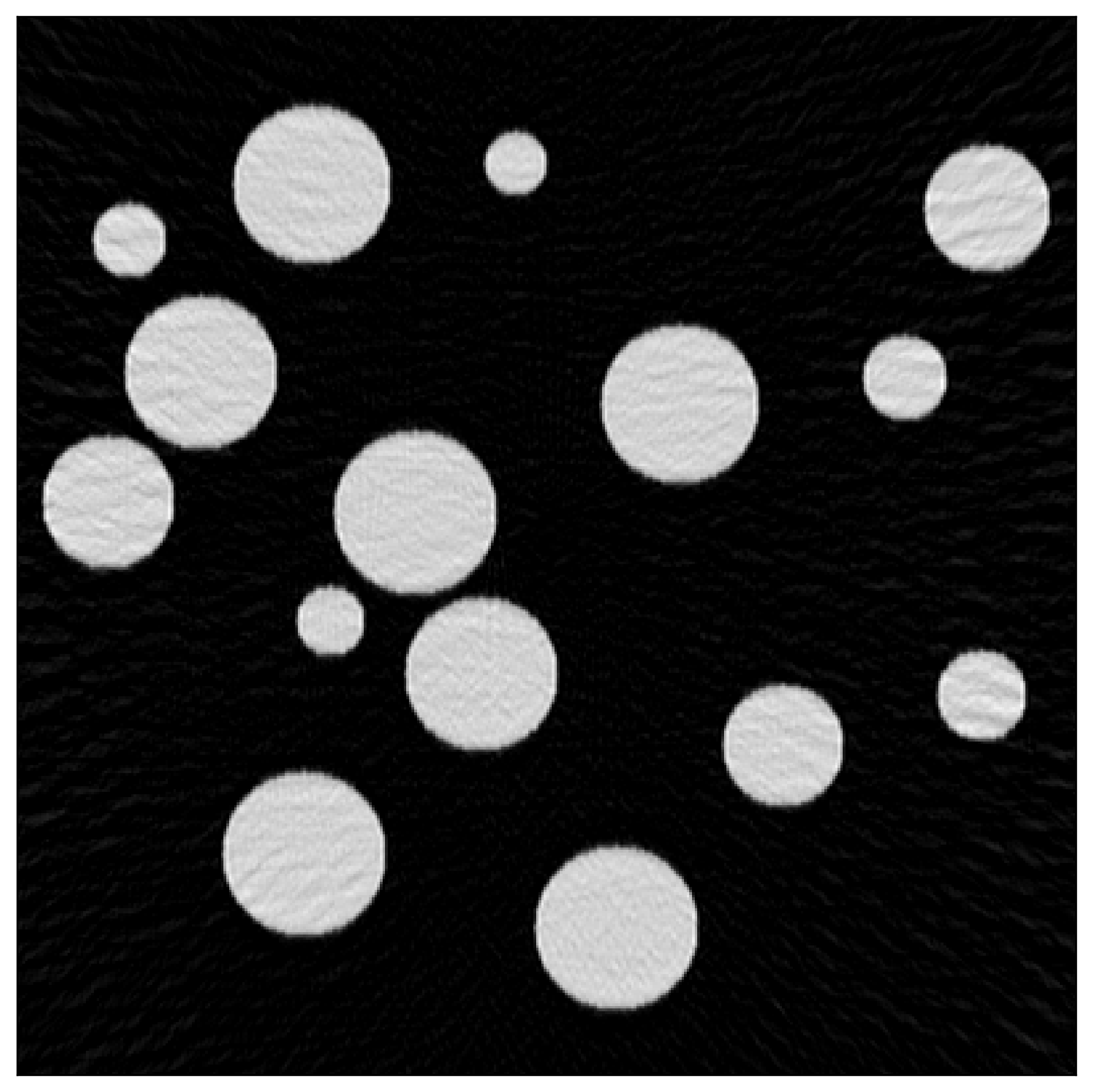}
        \caption*{Learned Method}
    \end{subfigure}

    \caption{Qualitative comparison showing reconstructions of a 2D phantom for the ground truth, FBP, and the learned method, respectively.}
    \label{fig:elliptical_eval}
\end{figure}

Table~\ref{tab:metrics_summary} compares the performance of the learned reconstruction method with the FBP algorithm, which is not inherently designed for elliptical trajectories. The quantitative metrics show that the learned filtered reconstruction significantly outperforms FBP in both MSE and SSIM, as the standard algorithm struggles to handle the non-static acquisition geometry. The qualitative results in figure \ref{fig:elliptical_eval} further highlight this improvement: while the standard FBP reconstruction (middle) exhibits significant artifacts, the learned approach (right) effectively compensates for the trajectory's geometry to produce a result nearly identical to the ground truth (left).

\begin{table}[t]
\centering
\caption{Reconstruction quality of the learned method against standard FBP for elliptical trajectory on the test set (mean $\pm$ std).}
\label{tab:metrics_summary}
\begin{tabular}{lcc}
\toprule
Method & MSE & SSIM \\
\midrule
Learned reconstruction & $0.0032 \pm 0.0002$ & $0.404 \pm 0.0127$ \\
FBP                   & $0.0535 \pm 0.0073$ & $0.122 \pm 0.0098$ \\
\bottomrule
\end{tabular}
\end{table}

\clearpage
\subsection{3D Circular Laminography - synthetic data}\label{Circular-lamino-section}
In the following, we optimize the filter and weights for filtered back-projection in a 3D circular laminography setup, illustrated in figure~\ref{fig:Lamino-setup}. This is achieved using the parameterization introduced in~(\ref{F_recon:circular}), incorporated into the optimization problem defined in~(\ref{supervised_impl}). The training data consist of synthetic 3D flat phantoms of size $40 \times 400 \times 400$ with three distinct layers. Each layer contains circles with randomly sampled positions and sizes. The acquisition comprises $100$ projections over a full $360^\circ$ rotation, each with a resolution of $512 \times 512$ pixels and the pixel size is set to $0.02$. The laminography angle is $\phi = 45^\circ$, with a SOD of $5$ and a SDD of $20$. 

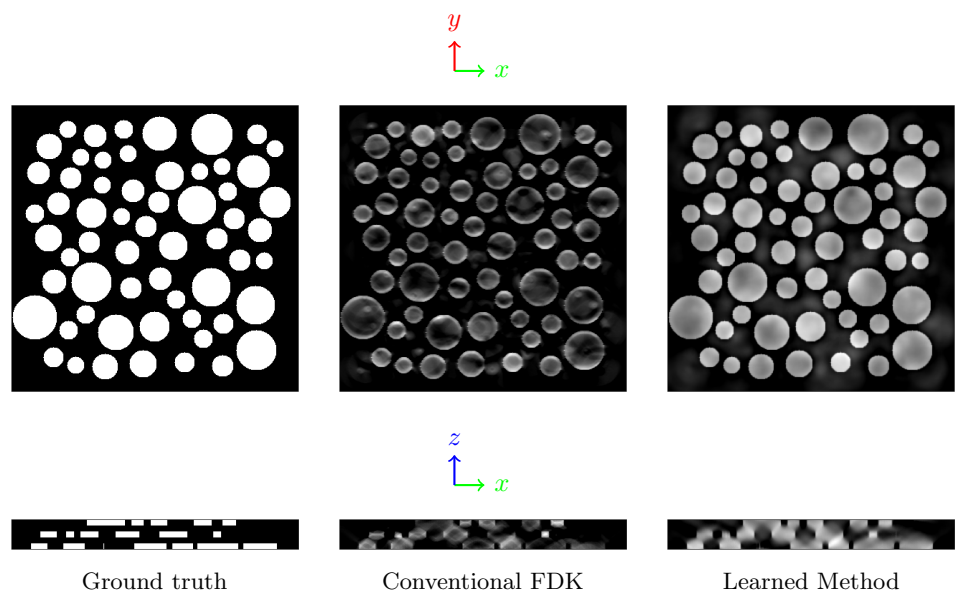
\begin{figure}[htbp]
    \centering
    


    \begin{tikzpicture}
      \begin{scope}[scale=0.4]
        \draw[->, thick, green] (0,0) -- (1,0) node[right] {$x$};
        \draw[->, thick, red]   (0,0) -- (0,1) node[above] {$y$};
      \end{scope}
    \end{tikzpicture}
    \vspace{0.5em}

    \begin{subfigure}{0.31\textwidth}
      \includegraphics[width=\linewidth]{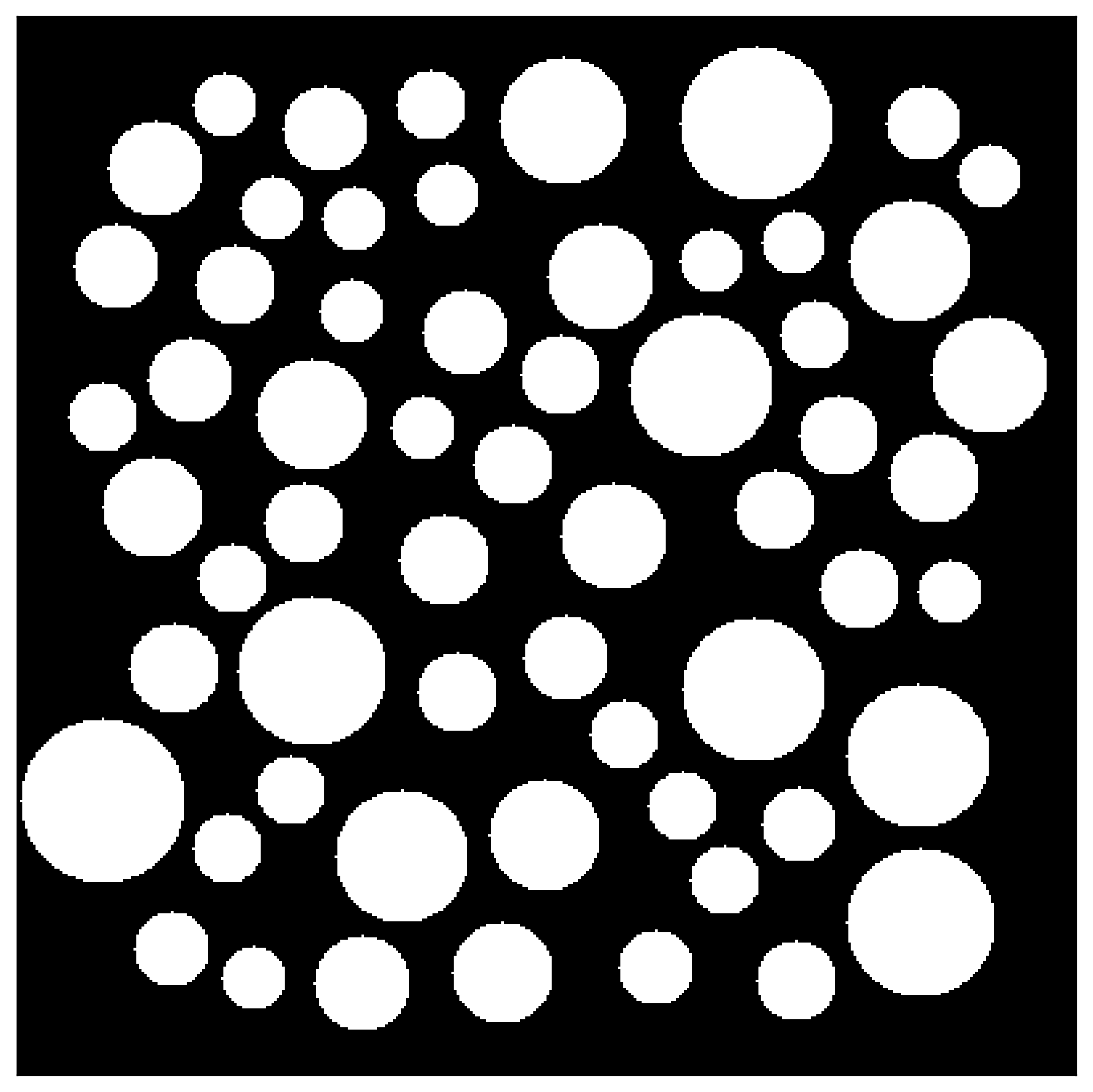}
    \end{subfigure}\hfill
    \begin{subfigure}{0.31\textwidth}
      \includegraphics[width=\linewidth]{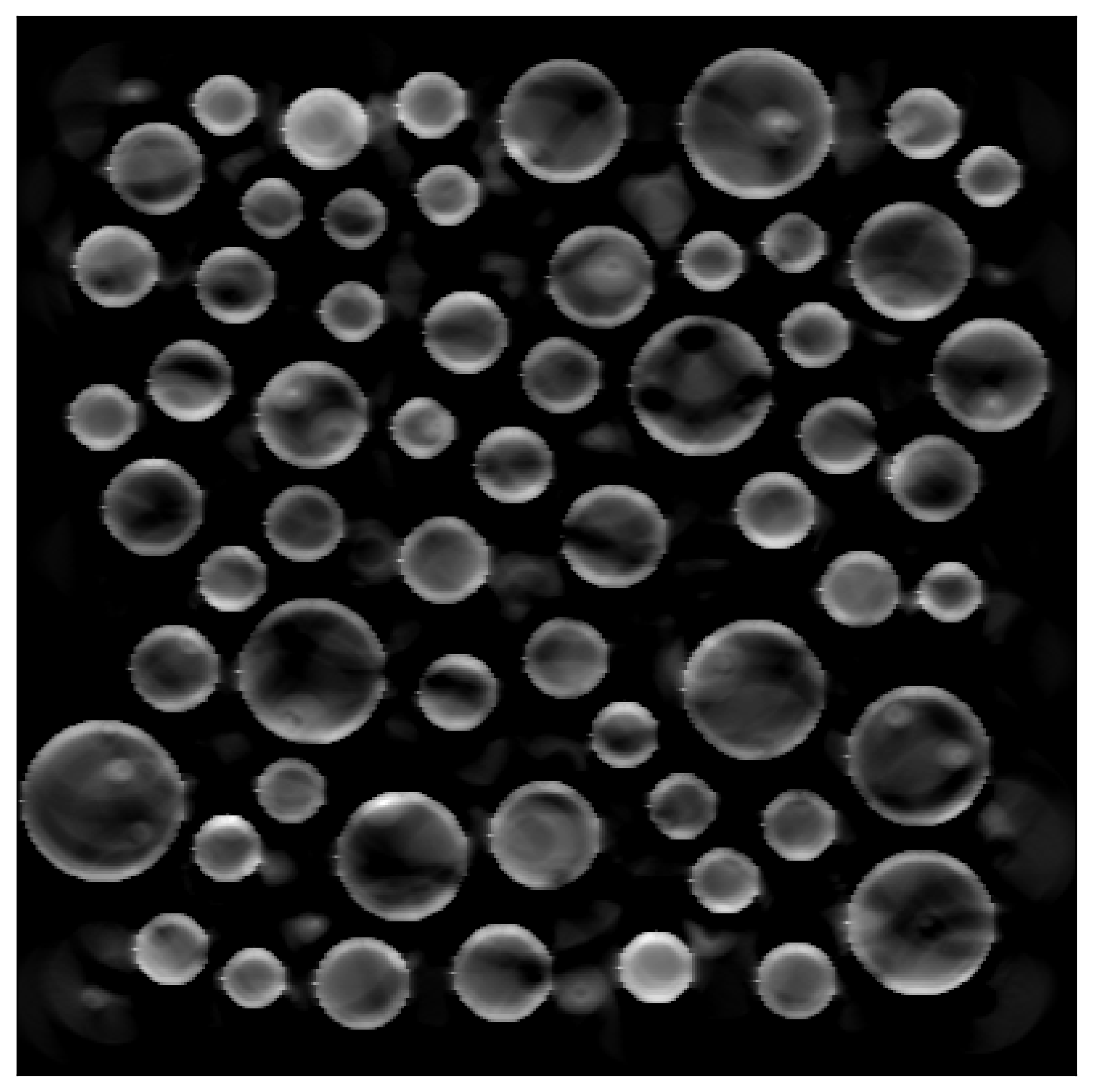}
    \end{subfigure}\hfill
    \begin{subfigure}{0.31\textwidth}
      \includegraphics[width=\linewidth]{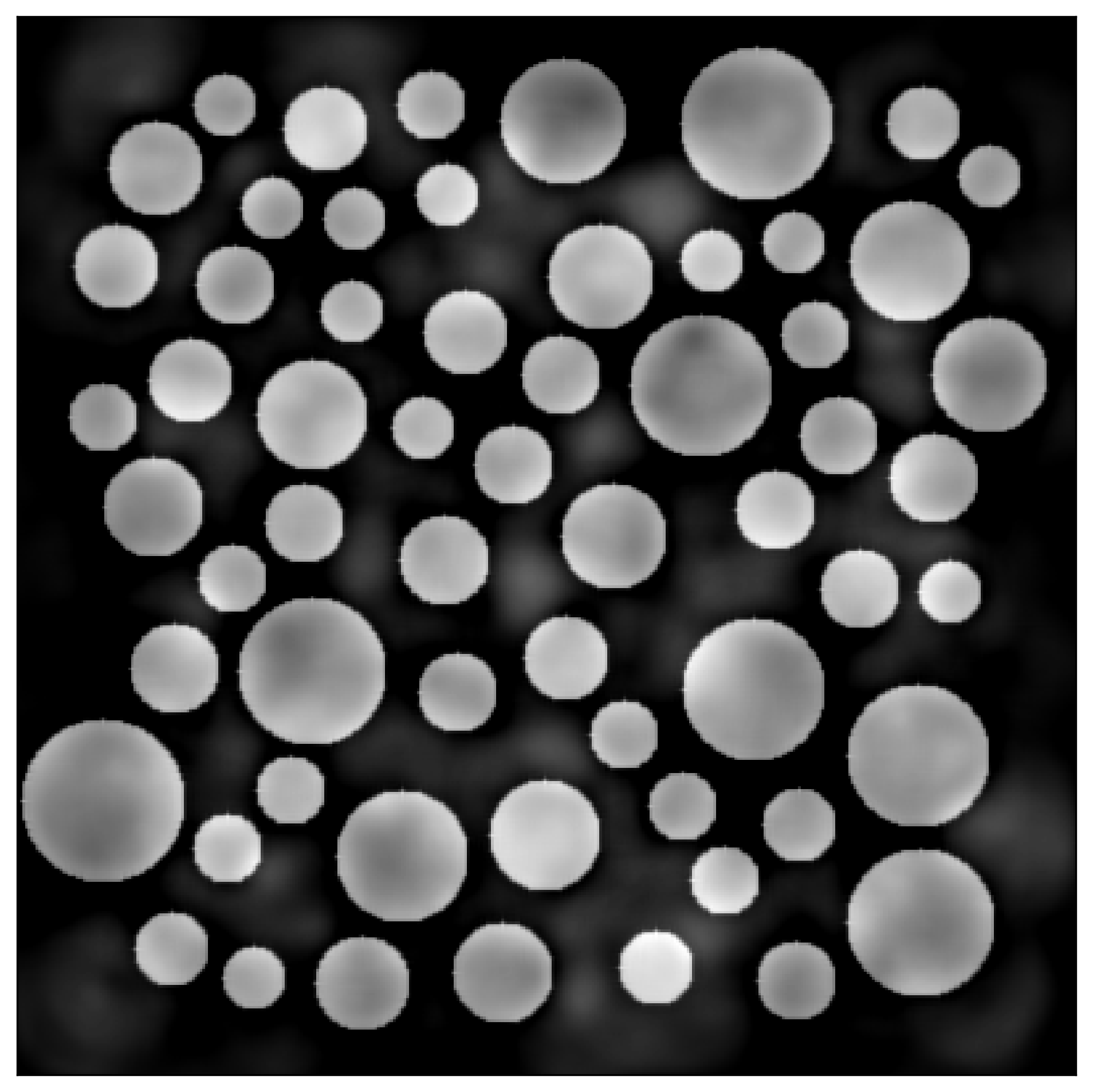}
    \end{subfigure}

    \vspace{1em}

    \begin{tikzpicture}
      \begin{scope}[scale=0.4]
        \draw[->, thick, green] (0,0) -- (1,0) node[right] {$x$};
        \draw[->, thick, blue]  (0,0) -- (0,1) node[above] {$z$};
      \end{scope}
    \end{tikzpicture}
    \vspace{0.5em}

    \begin{subfigure}{0.31\textwidth}
      \includegraphics[width=\linewidth]{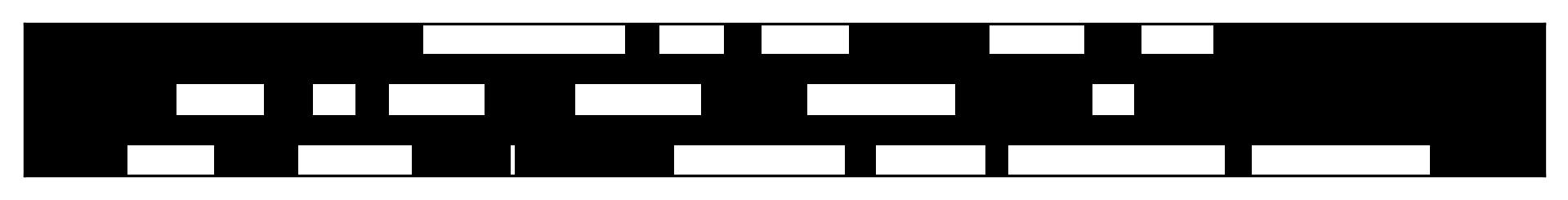}
      \caption*{Ground truth}
    \end{subfigure}\hfill
    \begin{subfigure}{0.31\textwidth}
      \includegraphics[width=\linewidth]{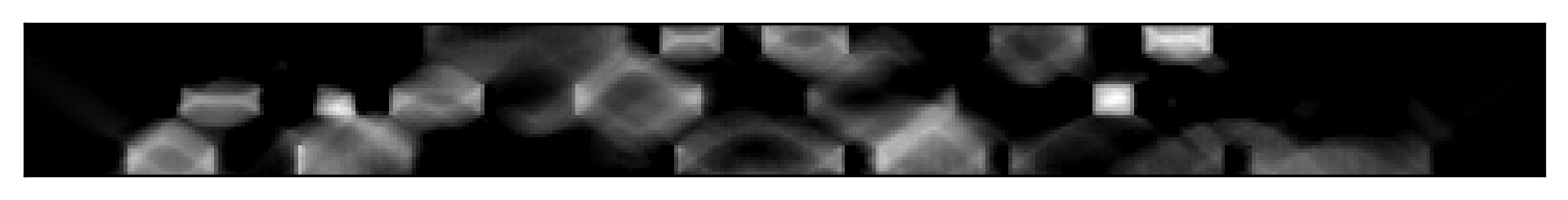}
    \caption*{Conventional FDK}
    \end{subfigure}\hfill
    \begin{subfigure}{0.31\textwidth}
      \includegraphics[width=\linewidth]{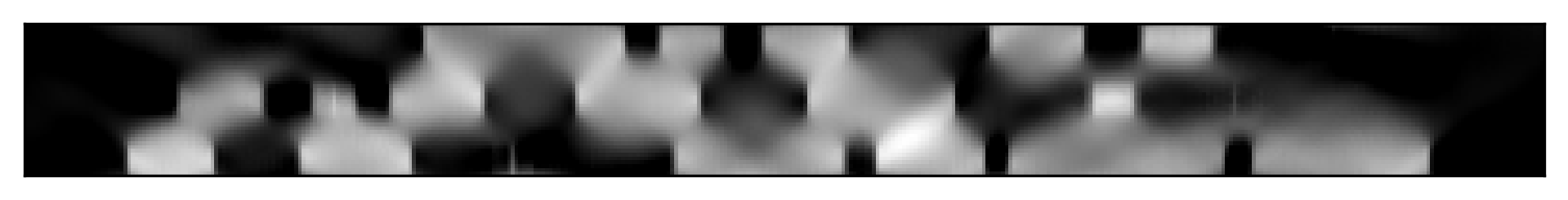}
    \caption*{Learned Method}
    \end{subfigure}

    \caption{Qualitative comparison of lateral and coronal slices in learned reconstruction (right column), respectively, against Ground Truth (left column) and conventional FDK (middle column).}
    \label{fig:lamino_eval}
\end{figure}

Table~\ref{tab:validation_metrics_} evaluates the reconstruction performance of the learned method compared to the conventional FDK algorithm quantitatively. It shows that the learned approach consistently outperforms FDK across the validation dataset, achieving significantly better MSE and SSIM distributions. These results are confirmed by the visual comparisons of a test phantom in figure~\ref{fig:lamino_eval}. While the conventional FDK reconstructions (middle column) exhibit severe artifacts in both lateral and coronal views, the learned reconstructions (right column) demonstrate substantial improvements in structural clarity, closely matching the ground truth (left column).

\begin{figure}[htbp]
    \centering
    
    \begin{subfigure}{0.45\textwidth}
        \centering
        \includegraphics[width=1\linewidth]{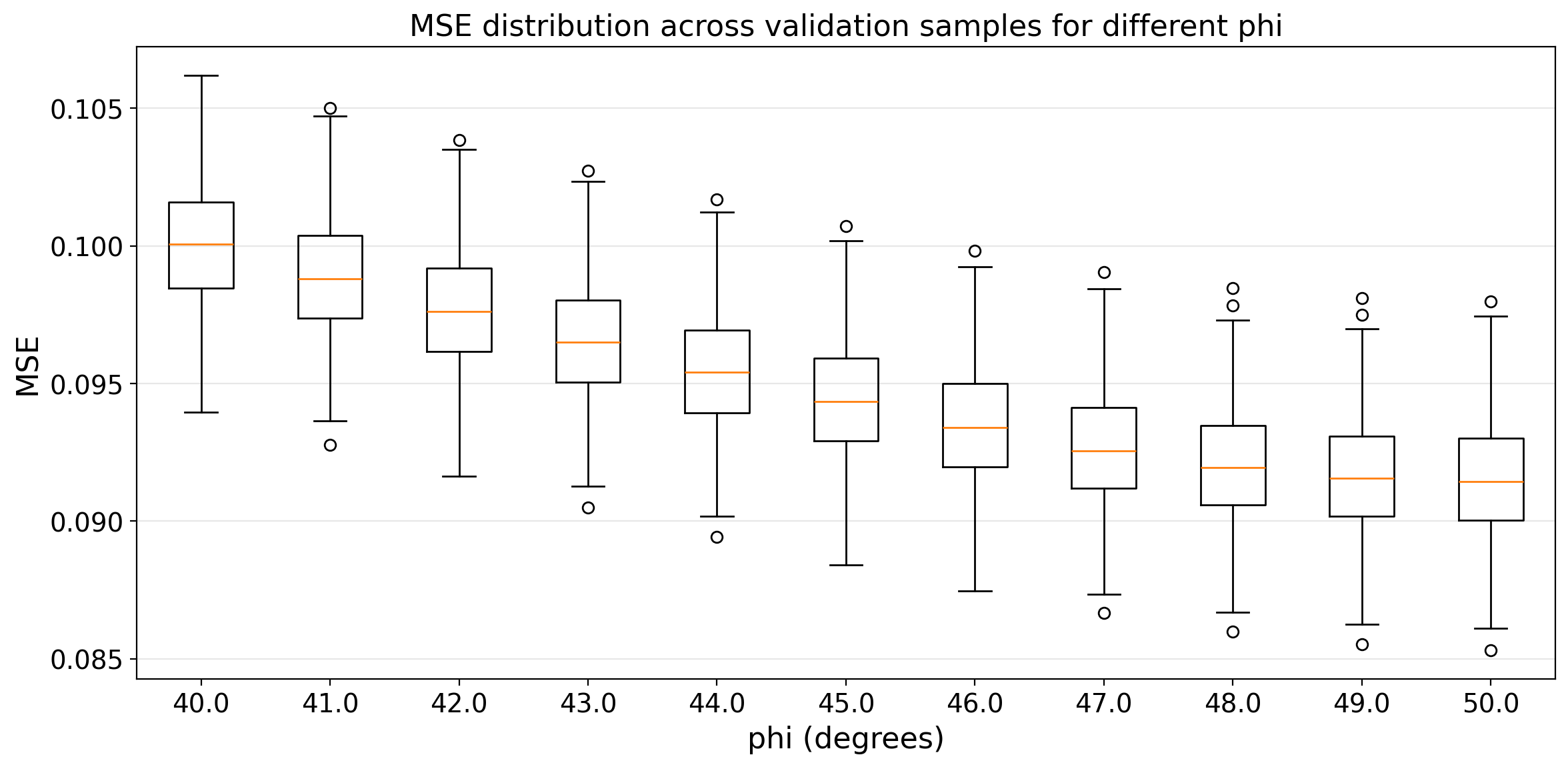}
        \caption*{MSE across angles}
    \end{subfigure}
    \hfill
    \begin{subfigure}{0.45\textwidth}
        \centering
        \includegraphics[width=1\linewidth]{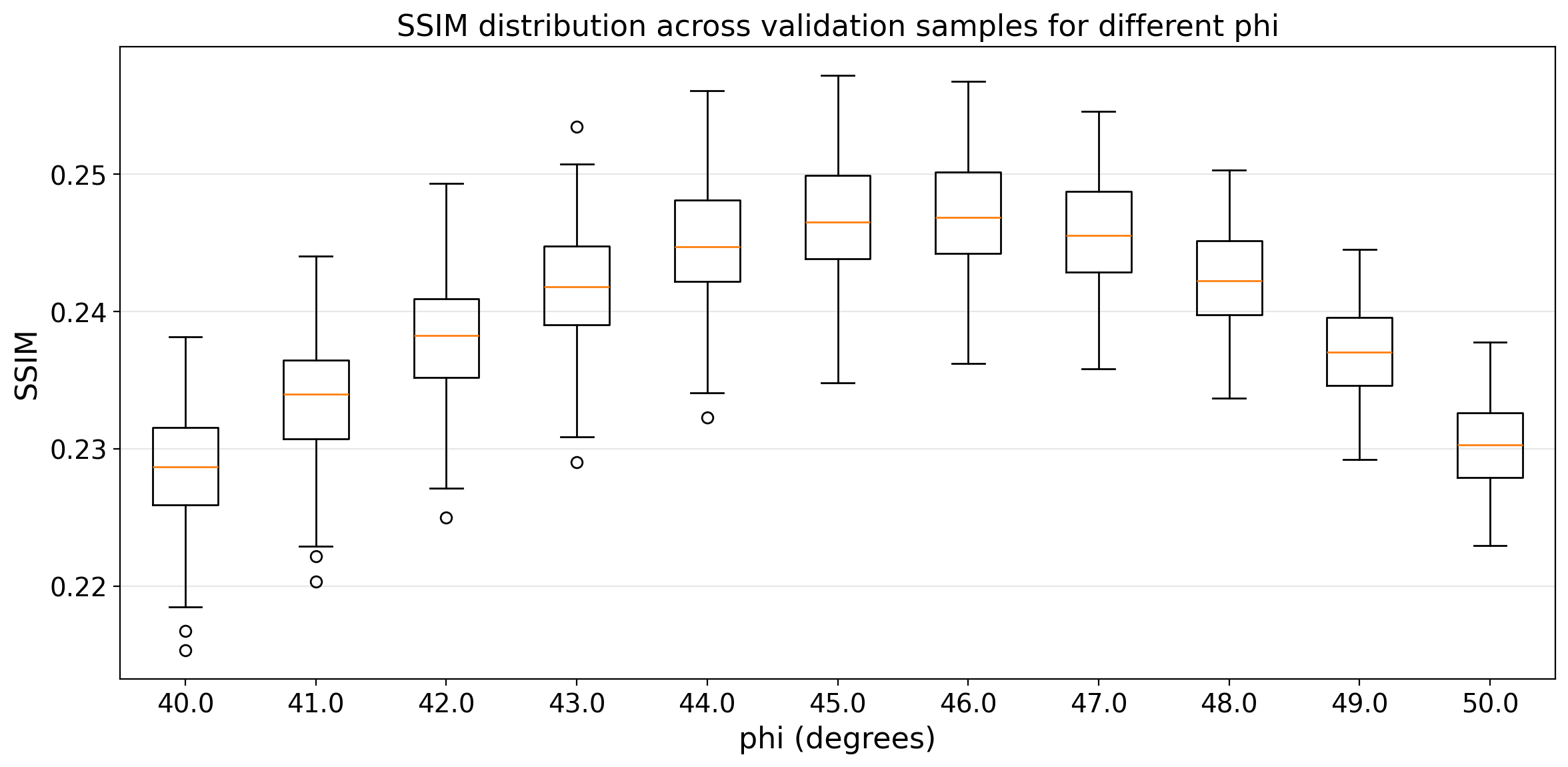}
        \caption*{SSIM across angles}
    \end{subfigure}

    \vspace{2em} 

    \begin{tikzpicture}[remember picture, baseline=(current bounding box.center)]
        \begin{scope}[scale=0.3]
            \draw[->, thick, green] (0,0,0) -- (1,0,0) node[right] {$x$};
            \draw[->, thick, red] (0,0,0) -- (0,1,0) node[above] {$y$};
        \end{scope}
    \end{tikzpicture}
    \begin{subfigure}{0.28\textwidth}
        \includegraphics[width=1\linewidth]{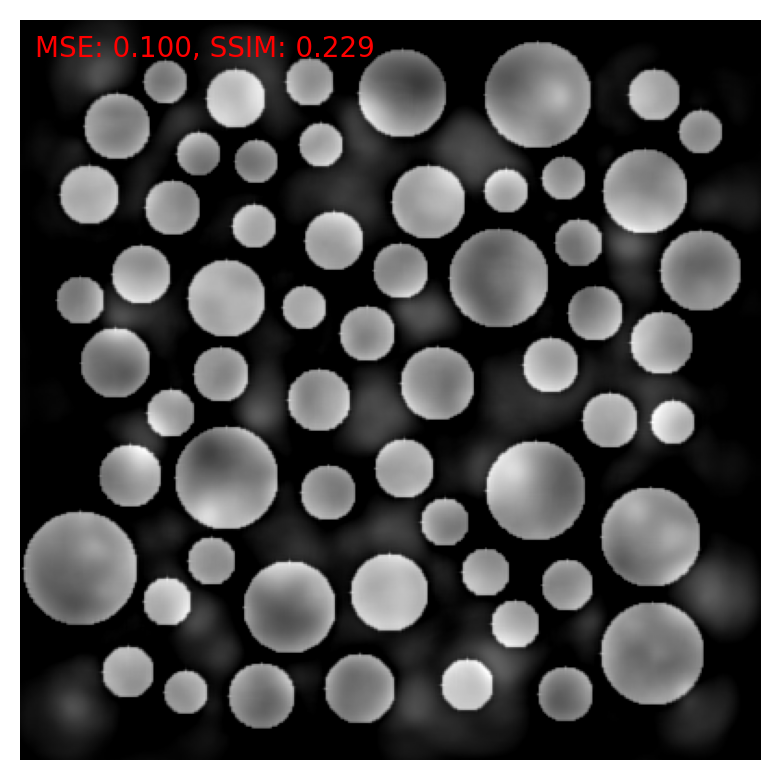}
    \end{subfigure}
    \begin{subfigure}{0.28\textwidth}
        \includegraphics[width=1\linewidth]{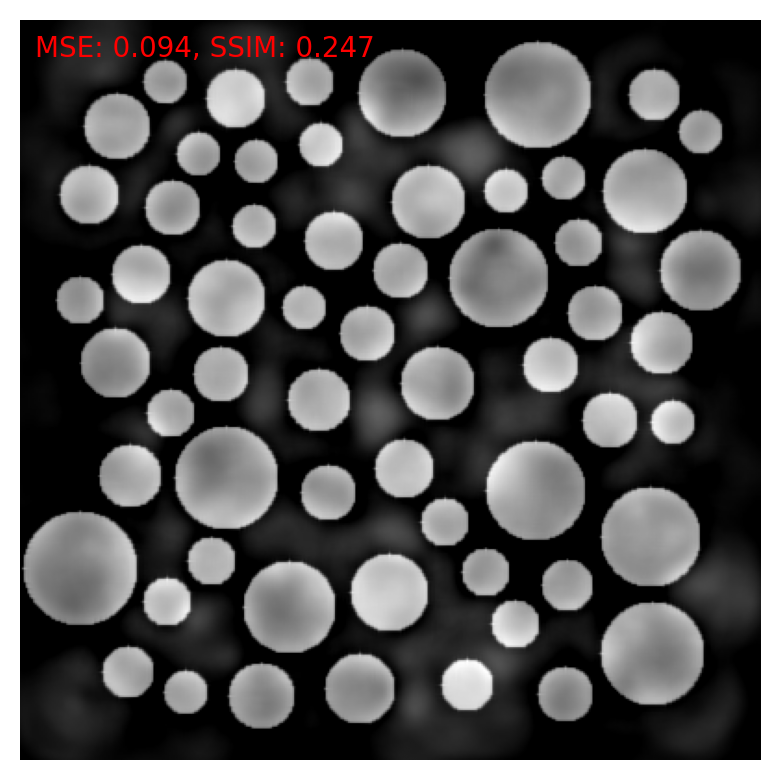}
    \end{subfigure}
    \begin{subfigure}{0.28\textwidth}
        \includegraphics[width=1\linewidth]{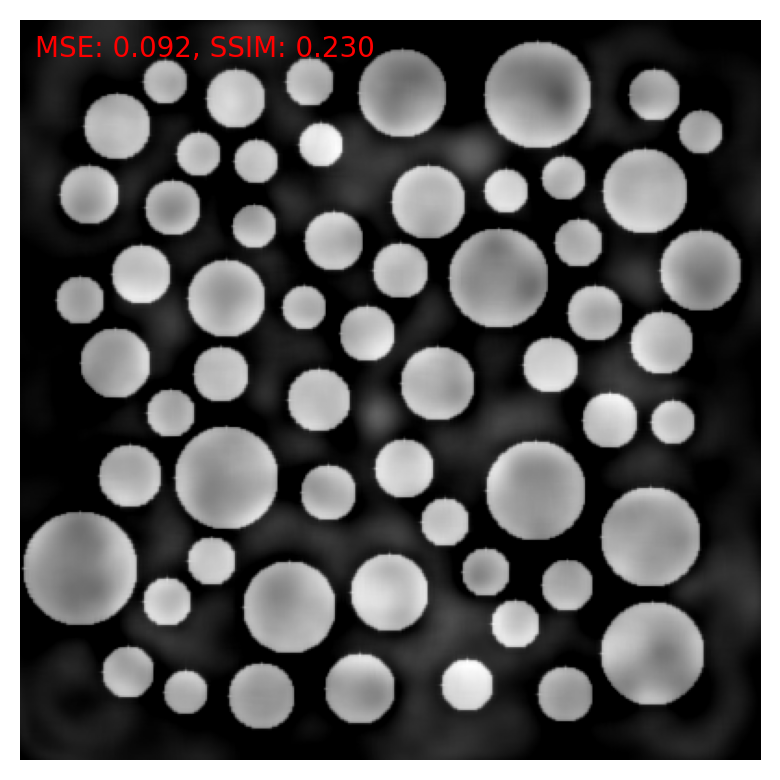}
    \end{subfigure}

    \vspace{1em} 

    \begin{tikzpicture}[remember picture, baseline=(current bounding box.center)]
        \begin{scope}[scale=0.3]
            \draw[->, thick, blue] (0,0,0) -- (0,1,0) node[above] {$z$};
            \draw[->, thick, green] (0,0,0) -- (1,0,0) node[right] {$x$};
        \end{scope}
    \end{tikzpicture}
    \begin{subfigure}{0.28\textwidth}
        \includegraphics[width=1\linewidth]{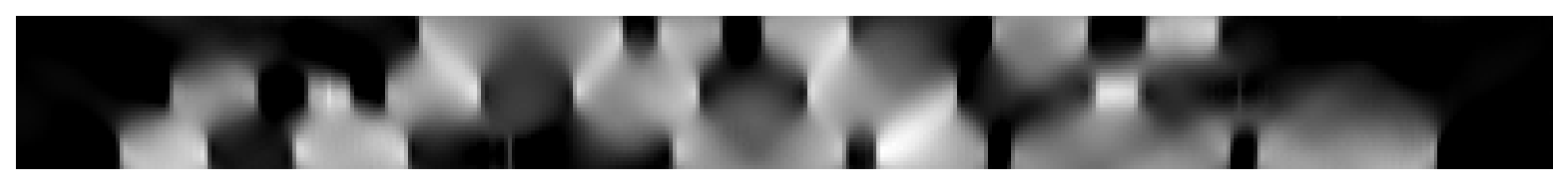}
        \caption*{$40^\circ$}        
    \end{subfigure}
    \begin{subfigure}{0.28\textwidth}
        \includegraphics[width=1\linewidth]{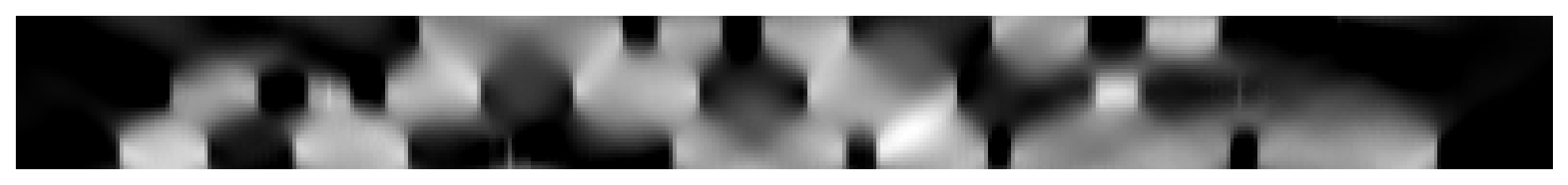}
        \caption*{$45^\circ$}        
    \end{subfigure}
    \begin{subfigure}{0.28\textwidth}
        \includegraphics[width=1\linewidth]{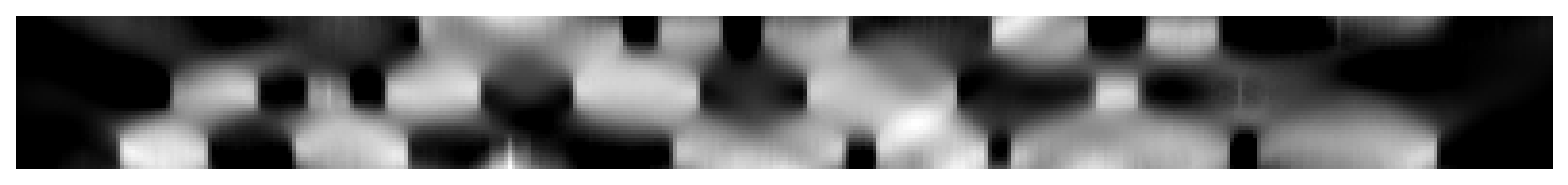}
        \caption*{$50^\circ$}    
    \end{subfigure}

  \caption{Generalizability of the learned parameters across varying laminography angles. \textbf{Top row:} Quantitative assessment via MSE and SSIM boxplots for variations up to $\pm5^\circ$. \textbf{Middle and bottom rows:} Qualitative lateral and coronal cross-sections, respectively, showing reconstructions at $40^\circ, 45^\circ$, and $50^\circ$.}
    \label{fig:angle_variation}
\end{figure}
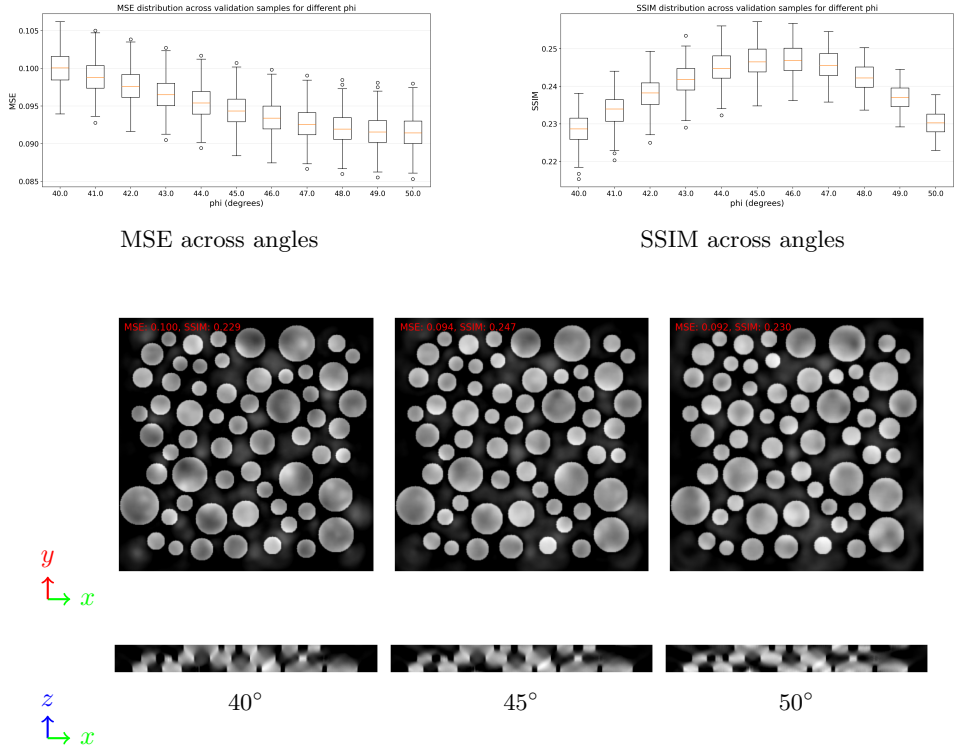

Figure~\ref{fig:angle_variation} evaluates the generalizability and robustness of the learned reconstruction parameters, which were trained at a $45^\circ$ laminography angle ($\phi$) and subsequently applied to measurements from different angles. The quantitative plots (top row) show that while the SSIM slightly declines as the acquisition angle deviates from the training setup, the MSE improves (decreases) at higher angles due to the enhanced depth information. The qualitative results (middle and bottom rows) visualize these effects through central slices in the lateral and coronal planes for $40^\circ, 45^\circ,$ and $50^\circ$, demonstrating stable performance across the tested range.

\begin{table}[t]
\centering
\caption{Evaluation of the learned method for cone-beam circular laminography against FDK results. It shows mean $\pm$ standard deviation.}
\label{tab:validation_metrics_}
\begin{tabular}{c|c|c}
\toprule
Method/Metric & MSE& SSIM \\
\midrule
Learned reconstruction & $0.094 \pm 0.002$ & $0.246 \pm 0.004$ \\
FDK   & $0.160 \pm 0.004$ & $0.059 \pm 0.003$ \\
\bottomrule
\end{tabular}
\end{table}

\clearpage
\subsection{3D Circular Laminography - real data}
To further validate the proposed method, we test it on an experimental cone-beam laminographic scan of a LEGO brick acquired by Waygate Technologies. The detector has a pixel size of $100\mu\text{m}$ and a native resolution of $3000 \times 3000$ pixels, downscaled by a factor of three to reduce computational cost. The acquisition comprises $72$ projections with imaging geometry defined by $SOD = 360$, $SDD = 1300$, and $\phi = 17.1^\circ$.

The model was first trained on the previously described synthetic dataset containing circle phantoms of size $125 \times 1000 \times 1000$ voxels while using experimental geometry. Then, learned filters and weights were applied to the LEGO brick dataset to evaluate generalization. For qualitative comparison, we present representative slices of the reconstructed volume alongside results obtained with conventional FDK and Nesterov-accelerated gradient descent (NAG), which is an iterative reconstruction method \cite{Nesterov}.

\begin{figure}[ht]
    \centering
    \begin{tikzpicture}[remember picture]
        \begin{scope}[scale=0.4]
            \draw[->, thick, green] (0,0,0) -- (1,0,0) node[right] {$x$};
            \draw[->, thick, red] (0,0,0) -- (0,1,0) node[above] {$y$};
        \end{scope}
    \end{tikzpicture}
    
    \makebox[0.3\textwidth]{FDK}%
    \makebox[0.3\textwidth]{Learned}%
    \makebox[0.3\textwidth]{Iterative}\\[0.3em]
    
    \begin{tikzpicture}[spy using outlines={rectangle, magnification=2, size=2.5cm, connect spies}]
        \node {\includegraphics[width=0.28\textwidth]{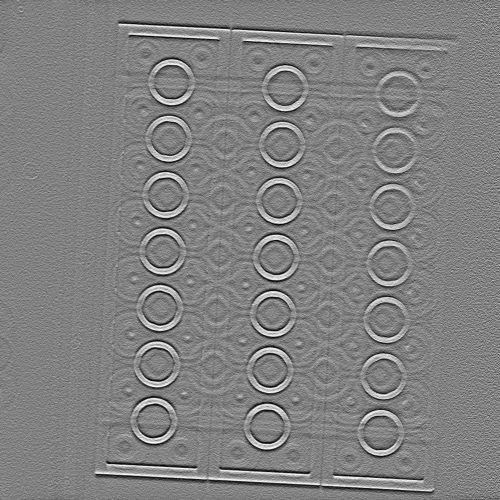}};
        \spy on (-1,1) in node [above] at (0,2.2);
    \end{tikzpicture}
    \hfill
    \begin{tikzpicture}[spy using outlines={rectangle, magnification=2, size=2.5cm, connect spies}]
        \node {\includegraphics[width=0.28\textwidth]{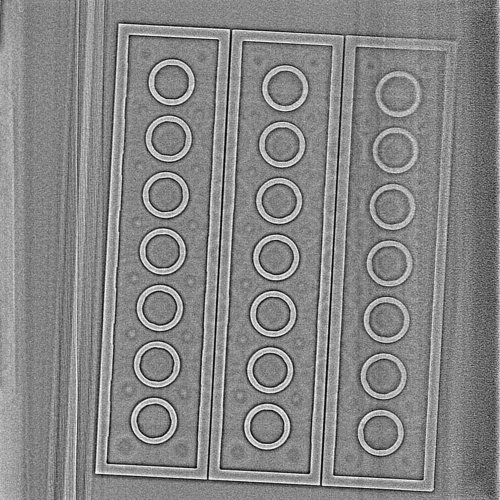}};
        \spy on (-1,1) in node [above] at (0,2.2);
    \end{tikzpicture}
    \hfill
    \begin{tikzpicture}[spy using outlines={rectangle, magnification=2, size=2.5cm, connect spies}]
        \node {\includegraphics[width=0.28\textwidth]{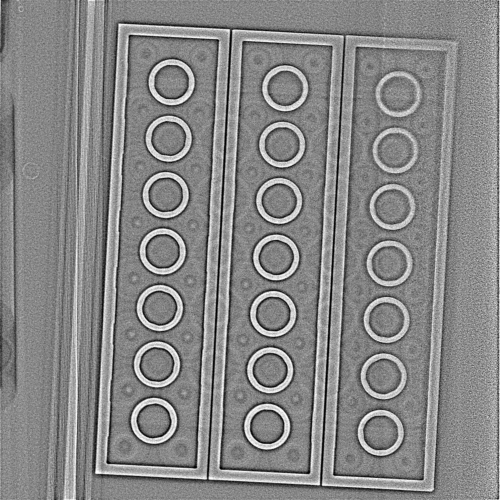}};
        \spy on (-1,1) in node [above] at (0,2.2);
    \end{tikzpicture}
    
    \caption{Comparison of the top axial slice through the reconstructions of experimental LEGO brick dataset. The coordinate system (top) applies to all images, with $x$ (red), $y$ (green), and $z$ (blue) axes. Each method shows the same ROI zoomed above the image.}
\label{top-layer:lamino}
\end{figure}
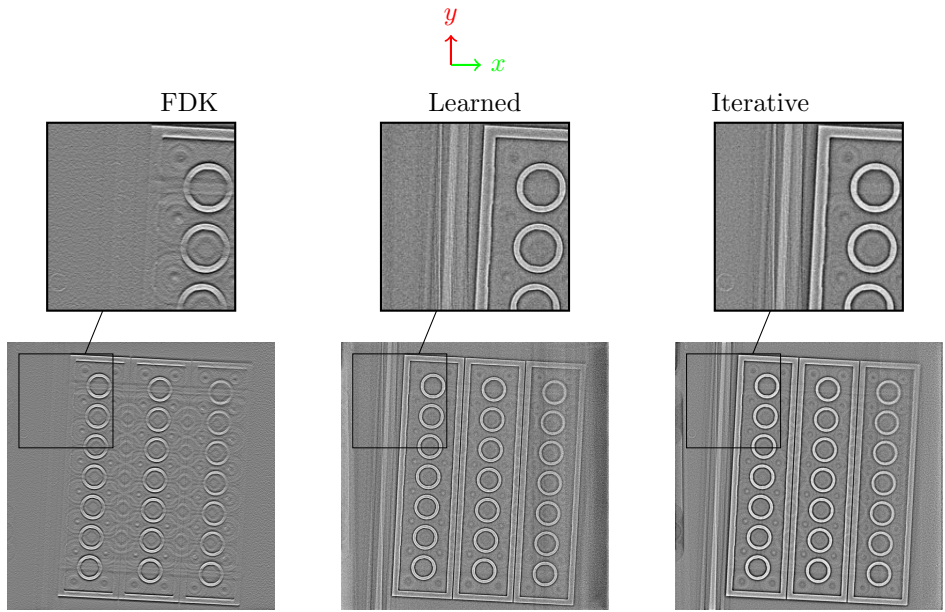

\begin{figure}
    \centering
    \begin{tikzpicture}[remember picture]
        \begin{scope}[scale=0.4]
            \draw[->, thick, green] (0,0,0) -- (1,0,0) node[right] {$x$};
            \draw[->, thick, red] (0,0,0) -- (0,1,0) node[above] {$y$};
        \end{scope}
    \end{tikzpicture}

    \makebox[0.3\textwidth]{FDK}%
    \makebox[0.3\textwidth]{Learned}%
    \makebox[0.3\textwidth]{Iterative}\\[0.3em]

    \begin{tikzpicture}[spy using outlines={rectangle, magnification=3, width=2.5cm, height=2.5cm, connect spies}]
        \node {\includegraphics[width=0.28\textwidth]{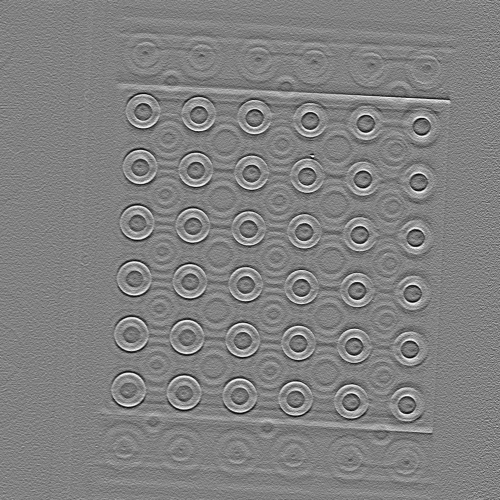}};
        \spy on (-1,0) in node [below] at (0,-2.2);
    \end{tikzpicture}
    \hfill
    \begin{tikzpicture}[spy using outlines={rectangle, magnification=3, width=2.5cm, height=2.5cm, connect spies}]
        \node {\includegraphics[width=0.28\textwidth]{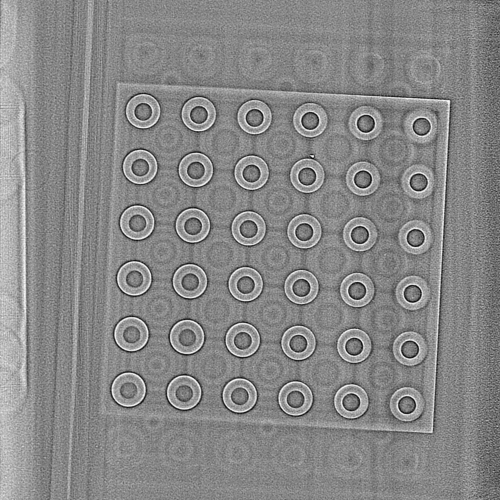}};
        \spy on (-1,0) in node [below] at (0,-2.2);
    \end{tikzpicture}
    \hfill
    \begin{tikzpicture}[spy using outlines={rectangle, magnification=3, width=2.5cm, height=2.5cm, connect spies}]
        \node {\includegraphics[width=0.28\textwidth]{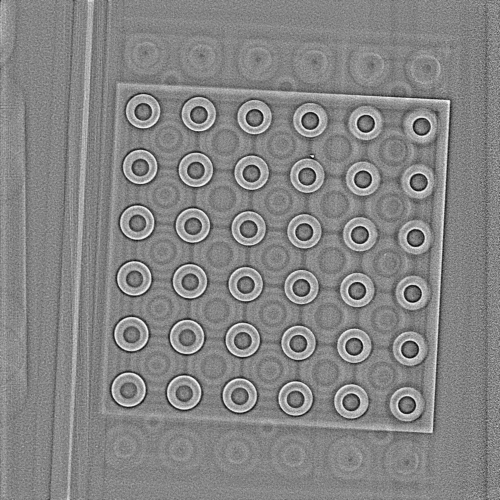}};
        \spy on (-1,0) in node [below] at (0,-2.2);
    \end{tikzpicture}
    \caption{Comparison of the bottom axial slice through the reconstructions of experimental LEGO brick dataset. Each method (FDK, learned, iterative) shows the same ROI zoomed above the image.}
    \label{bottom-layer:lamino}
\end{figure}

\begin{figure}
    \centering
    \begin{tikzpicture}[remember picture]
        \begin{scope}[scale=0.4]
            \draw[->, thick, blue] (0,0,0) -- (0,1,0) node[above] {$z$};
            \draw[->, thick, green] (0,0,0) -- (1,0,0) node[right] {$x$};
        \end{scope}
    \end{tikzpicture}

    \makebox[0.3\textwidth]{FDK}%
    \makebox[0.3\textwidth]{Learned}%
    \makebox[0.3\textwidth]{Iterative}\\[0.3em]
    
    \begin{tikzpicture}[spy using outlines={rectangle, magnification=3, width=2.5cm, height=1cm, connect spies}]
        \node {\includegraphics[width=0.28\textwidth]{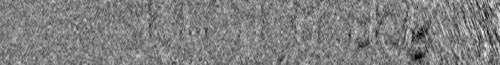}};
        \spy on (0,0) in node [above] at (0,0.5);
    \end{tikzpicture}
    \hfill
    \begin{tikzpicture}[spy using outlines={rectangle, magnification=3, width=2.5cm, height=1cm, connect spies}]
        \node {\includegraphics[width=0.28\textwidth]{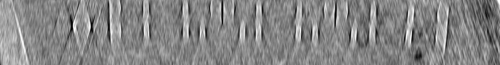}};
        \spy on (0,0) in node [above] at (0,0.5);
    \end{tikzpicture}
    \hfill
    \begin{tikzpicture}[spy using outlines={rectangle, magnification=3, width=2.5cm, height=1cm, connect spies}]
        \node {\includegraphics[width=0.28\textwidth]{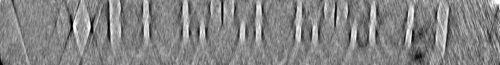}};
        \spy on (0,0) in node [above] at (0,0.5);
    \end{tikzpicture}
    \caption{Comparison of the central coronal slice through the reconstructions of experimental LEGO brick dataset. Each method (FDK, learned, iterative) shows the same ROI zoomed above the image.}
    \label{coronal-view:lamino}
\end{figure}

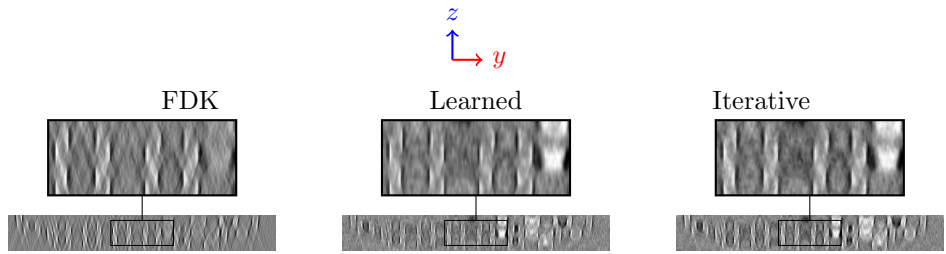
\begin{figure}
    \centering
    \begin{tikzpicture}[remember picture]
        \begin{scope}[scale=0.4]
            \draw[->, thick, blue] (0,0,0) -- (0,1,0) node[above] {$z$};
            \draw[->, thick, red] (0,0,0) -- (1,0,0) node[right] {$y$};
        \end{scope}
    \end{tikzpicture}
    
    \makebox[0.3\textwidth]{FDK}%
    \makebox[0.3\textwidth]{Learned}%
    \makebox[0.3\textwidth]{Iterative}\\[0.3em]
    \begin{tikzpicture}[spy using outlines={rectangle, magnification=3, width=2.5cm, height=1cm, connect spies}]
        \node {\includegraphics[width=0.28\textwidth]{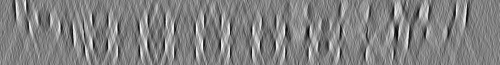}};
        \spy on (0,0) in node [above] at (0,0.5);
    \end{tikzpicture}
    \hfill
    \begin{tikzpicture}[spy using outlines={rectangle, magnification=3, width=2.5cm, height=1cm, connect spies}]
        \node {\includegraphics[width=0.28\textwidth]{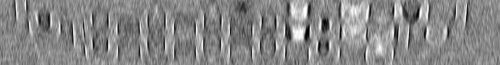}};
        \spy on (0,0) in node [above] at (0,0.5);
    \end{tikzpicture}
    \hfill
    \begin{tikzpicture}[spy using outlines={rectangle, magnification=3, width=2.5cm, height=1cm, connect spies}]
        \node {\includegraphics[width=0.28\textwidth]{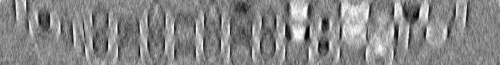}};
        \spy on (0,0) in node [above] at (0,0.5);
    \end{tikzpicture}
    \caption{Comparison of the central sagittal slice through the reconstructions of experimental LEGO brick dataset. Each method (FDK, learned, iterative) shows the same ROI zoomed above the image.}
    \label{sagittal-view:lamino}
\end{figure}

Figure~\ref{top-layer:lamino} shows axial views of the top layer reconstructed with the three methods. Treating the NAG results as a reference, it can be observed that conventional FDK fails to recover several structures, whereas the learned reconstruction preserves those features, closely matching the NAG solution. A similar observation is made for the bottom layer, shown in figure~\ref{bottom-layer:lamino}. Coronal views of the central slice are presented in figure~\ref{coronal-view:lamino}, where FDK performs poorly compared to both the learned and iterative reconstructions. Finally, the sagittal view in figure~\ref{sagittal-view:lamino} indicates that although FDK appears more comparable in this orientation, the learned reconstruction still demonstrates clear suppression of streak artifacts and achieves image quality on par with NAG.

Importantly, this result is obtained without iterative refinement at the reconstruction step, resulting in a comparable computational cost to the conventional FDK. The learned reconstruction not only outperforms analytical FDK on simulated phantoms but also generalizes effectively to experimental data, producing reconstructions comparable to iterative algorithms while retaining the computational efficiency of a direct method.

\clearpage
\section{Conclusion}
This paper introduces a data-driven framework to optimize filters and weighting factors for the filtered back-projection method, specifically targeting the limitations of traditional analytical techniques in challenging imaging conditions, such as with non-conventional geometries, significant noise and undersampling. By formulating a regularized optimization problem, the approach learns noise and geometry-adapted parameters directly from training data to address the ill-posedness of the reconstruction. Unlike the fixed filters used in conventional FBP and FDK, the learned filters act as spectral gains that adaptively attenuate high frequencies to suppress noise-induced artifacts while preserving structural details. The proposed method demonstrates significant performance gains in both 2D and 3D scenarios, including parallel-beam, fan-beam, and complex non-circular trajectories like elliptical paths and circular laminography. Quantitative results show that these learned angle-dependent filters and weights consistently lead to lower reconstruction errors and higher structural similarity compared to classical FBP and FDK methods. The framework enables the effective handling of incomplete sampling and acquisition geometries that deviate from standard CT configurations, which typically lead to strong artifacts in analytical reconstructions. Furthermore, the approach proves its practical utility by successfully generalizing from synthetic training data to real-world experimental measurements, such as the LEGO brick dataset. Ultimately, this framework offers a flexible and robust alternative that retains the real-time efficiency of direct back-projection while achieving the superior image quality associated with advanced iterative methods.
\section*{Acknowledgements}
The authors thank Dr.~Alexander Suppes (Waygate Technologies) for acquiring the LEGO brick laminography dataset. This work was supported by the European Union’s Horizon 2020 research and innovation programme under the Marie Sk\l{}odowska--Curie
grant agreement No.~956172 (xCTing).
\clearpage
\appendix

\section{Stability}
\label{Stability_proof}
\begin{proof}
Let $\pi$ and $\pi'$ be two joint distributions of $(\mathbf{x},\mathbf{y})$ satisfying Assumptions (A1)--(A3), and set
\[
\mathbf{\Sigma}_{\mathbf{yy}}(\pi)=\mathbb{E}_{\pi}[\mathbf{y}\mathbf{y}^\top],\qquad
\mathbf{\Sigma}_{\mathbf{xy}}(\pi)=\mathbb{E}_{\pi}[\mathbf{x}\mathbf{y}^\top],
\]
with analogous definitions for $\pi'$. Denote the objective in~\eqref{Learned Inverse_supervised} by $F_\pi(\mathbf{B})$.
Its gradient reads
\begin{equation}\label{GradF_stability_clean}
\nabla F_\pi(\mathbf{B})
=
2\,\mathbf{A}\mathbf{A}^\top \mathbf{B}\,\mathbf{\Sigma}_{\mathbf{yy}}(\pi)
-2\,\mathbf{A}\,\mathbf{\Sigma}_{\mathbf{xy}}(\pi)
+\lambda \nabla\rho(\mathbf{B}).
\end{equation}
Let $\mathbf{B}^\ast_\pi$ and $\mathbf{B}^\ast_{\pi'}$ be minimizers of $F_\pi$ and $F_{\pi'}$, hence
$\nabla F_\pi(\mathbf{B}^\ast_\pi)=\nabla F_{\pi'}(\mathbf{B}^\ast_{\pi'})=\mathbf{0}$.

\smallskip
\noindent\textbf{Step 1 (strong convexity).}
Under (A1)--(A2) the data term is strongly convex in $\mathbf{B}$, and since $\rho$ is $\mu$-strongly convex,
$F_\pi$ is $(\gamma+\lambda\mu)$-strongly convex for some $\gamma>0$ \cite{nesterov2013introductory}. Therefore,
for any $\mathbf{B},\mathbf{B}'$,
\[
(\gamma+\lambda\mu)\,\|\mathbf{B}-\mathbf{B}'\|_F^2
\le
\big\langle \nabla F_\pi(\mathbf{B})-\nabla F_\pi(\mathbf{B}'),\, \mathbf{B}-\mathbf{B}'\big\rangle .
\]
With $\mathbf{B}=\mathbf{B}^\ast_\pi$ and $\mathbf{B}'=\mathbf{B}^\ast_{\pi'}$ and using
$\nabla F_\pi(\mathbf{B}^\ast_\pi)=\mathbf{0}$, we obtain
\begin{equation}\label{SC_formula_stability_clean}
(\gamma+\lambda\mu)\,\|\mathbf{B}^\ast_\pi-\mathbf{B}^\ast_{\pi'}\|_F^2
\le
-\big\langle \nabla F_\pi(\mathbf{B}^\ast_{\pi'}),\, \mathbf{B}^\ast_\pi-\mathbf{B}^\ast_{\pi'}\big\rangle.
\end{equation}

\smallskip
\noindent\textbf{Step 2 (difference of gradients).}
Using~\eqref{GradF_stability_clean} and $\nabla F_{\pi'}(\mathbf{B}^\ast_{\pi'})=\mathbf{0}$ gives
\[
\nabla F_\pi(\mathbf{B}^\ast_{\pi'})
=
2\,\mathbf{A}\mathbf{A}^\top \mathbf{B}^\ast_{\pi'}\!\left(\mathbf{\Sigma}_{\mathbf{yy}}(\pi)-\mathbf{\Sigma}_{\mathbf{yy}}(\pi')\right)
-2\,\mathbf{A}\!\left(\mathbf{\Sigma}_{\mathbf{xy}}(\pi)-\mathbf{\Sigma}_{\mathbf{xy}}(\pi')\right).
\]
Hence,
\begin{equation}\label{grad_bound_clean}
\|\nabla F_\pi(\mathbf{B}^\ast_{\pi'})\|_F
\le
2\|\mathbf{A}\mathbf{A}^\top\|_F\,\|\mathbf{B}^\ast_{\pi'}\|_F\,\|\mathbf{\Sigma}_{\mathbf{yy}}(\pi)-\mathbf{\Sigma}_{\mathbf{yy}}(\pi')\|_F
+
2\|\mathbf{A}\|_F\,\|\mathbf{\Sigma}_{\mathbf{xy}}(\pi)-\mathbf{\Sigma}_{\mathbf{xy}}(\pi')\|_F .
\end{equation}

\smallskip
\noindent\textbf{Step 3 (moment bounds via $W_2$).}
Using the coupling characterization of $W_2$ \cite{villani2008optimal} together with the uniform moment bounds in (A1),
there exists $L>0$ such that
\[
\|\mathbf{\Sigma}_{\mathbf{yy}}(\pi)-\mathbf{\Sigma}_{\mathbf{yy}}(\pi')\|_F \le L\,W_2(\pi,\pi'),
\qquad
\|\mathbf{\Sigma}_{\mathbf{xy}}(\pi)-\mathbf{\Sigma}_{\mathbf{xy}}(\pi')\|_F \le L\,W_2(\pi,\pi').
\]
Moreover, by strong convexity and the uniform moment bounds, $\|\mathbf{B}^\ast_{\pi'}\|_F$ is uniformly bounded by some
$M>0$. Substituting these bounds into~\eqref{grad_bound_clean} yields
\[
\|\nabla F_\pi(\mathbf{B}^\ast_{\pi'})\|_F
\le
2L\big(\|\mathbf{A}\mathbf{A}^\top\|_F\,M+\|\mathbf{A}\|_F\big)\,W_2(\pi,\pi')
=:C_1\,W_2(\pi,\pi').
\]

\smallskip
\noindent\textbf{Step 4 (conclusion).}
From~\eqref{SC_formula_stability_clean} and Cauchy--Schwarz,
\[
(\gamma+\lambda\mu)\,\|\mathbf{B}^\ast_\pi-\mathbf{B}^\ast_{\pi'}\|_F^2
\le
\|\nabla F_\pi(\mathbf{B}^\ast_{\pi'})\|_F\,\|\mathbf{B}^\ast_\pi-\mathbf{B}^\ast_{\pi'}\|_F,
\]
and therefore
\[
\|\mathbf{B}^\ast_\pi-\mathbf{B}^\ast_{\pi'}\|_F
\le
\frac{\|\nabla F_\pi(\mathbf{B}^\ast_{\pi'})\|_F}{\gamma+\lambda\mu}
\le
\frac{C_1}{\gamma+\lambda\mu}\,W_2(\pi,\pi').
\]
Since $\gamma>0$ is fixed (for given $\mathbf{A}$ and the moment bounds), we may absorb it into the constant and write
\[
\|\mathbf{B}^\ast_\pi-\mathbf{B}^\ast_{\pi'}\|_F
\le
\frac{C}{\lambda\mu}\,W_2(\pi,\pi'),
\]
for some $C>0$ independent of $\lambda$.
\end{proof}

\section{Reconstruction error bound}\label{Reconstruction_error}
\begin{proof}
Let $\mathbf{y}=\mathbf{A}\mathbf{x}+\mathbf{n}$, where the noise satisfies
$\mathbb{E}[\mathbf{n}\mid \mathbf{x}]=\mathbf{0}$.
Assume that $\mathbf{A}$ has full row rank (\textbf{A2}), so that the Moore--Penrose pseudoinverse $\mathbf{A}^\dagger=\mathbf{A}^\top(\mathbf{A}\mathbf{A}^\top)^{-1}$ is well defined and that
$\mathbf{A}^\dagger\mathbf{A}$ is the orthogonal projector onto the range of
$\mathbf{A}^\top$.

The reconstruction error reads
\[
\mathbf{A}^\top \mathbf{B}\mathbf{y}-\mathbf{x}
=
(\mathbf{A}^\top \mathbf{B}\mathbf{A}-\mathbf{I})\mathbf{x}
+
\mathbf{A}^\top \mathbf{B}\mathbf{n}.
\]
Adding and subtracting $\mathbf{A}^\dagger\mathbf{A}\mathbf{x}$ yields
\[
\mathbf{A}^\top \mathbf{B}\mathbf{y}-\mathbf{x}
=
\mathbf{A}^\top \mathbf{B}\mathbf{n}
+
(\mathbf{A}^\top\mathbf{B}-\mathbf{A}^\dagger)\mathbf{A}\mathbf{x}
+
(\mathbf{A}^\dagger\mathbf{A}-\mathbf{I})\mathbf{x}.
\]

Since both $\mathbf{A}^\top \mathbf{B}\mathbf{n}$ and
$(\mathbf{A}^\top\mathbf{B}-\mathbf{A}^\dagger)\mathbf{A}\mathbf{x}$
belong to the range of $\mathbf{A}^\top$, while
$(\mathbf{A}^\dagger\mathbf{A}-\mathbf{I})\mathbf{x}
=-(\mathbf{I}-\mathbf{A}^\dagger\mathbf{A})\mathbf{x}$ belongs to the null-space of $\mathbf{A}$,
then the third term is orthogonal to the sum of first two and the Pythagorean theorem gives
\[
\|\mathbf{A}^\top \mathbf{B}\mathbf{y}-\mathbf{x}\|_2^2
=
\|\mathbf{A}^\top \mathbf{B}\mathbf{n}
+
(\mathbf{A}^\top\mathbf{B}-\mathbf{A}^\dagger)\mathbf{A}\mathbf{x}\|_2^2
+
\|(\mathbf{A}^\dagger\mathbf{A}-\mathbf{I})\mathbf{x}\|_2^2 .
\]

Expanding the first term and taking expectation,
\begin{align*}
\mathbb{E}\|\mathbf{A}^\top \mathbf{B}\mathbf{y}-\mathbf{x}\|_2^2
&=
\mathbb{E}\|\mathbf{A}^\top \mathbf{B}\mathbf{n}\|_2^2
+
\mathbb{E}\|(\mathbf{A}^\top\mathbf{B}-\mathbf{A}^\dagger)\mathbf{A}\mathbf{x}\|_2^2 \\
&\qquad
+
2\,\mathbb{E}\!\left[
\left\langle
\mathbf{A}^\top \mathbf{B}\mathbf{n},
(\mathbf{A}^\top\mathbf{B}-\mathbf{A}^\dagger)\mathbf{A}\mathbf{x}
\right\rangle
\right]
+
\mathbb{E}\|(\mathbf{A}^\dagger\mathbf{A}-\mathbf{I})\mathbf{x}\|_2^2 .
\end{align*}

Using conditional expectation and $\mathbb{E}[\mathbf{n}\mid\mathbf{x}]=0$,
the cross term vanishes:
\[
\mathbb{E}\!\left[
\left\langle
\mathbf{A}^\top \mathbf{B}\mathbf{n},
(\mathbf{A}^\top\mathbf{B}-\mathbf{A}^\dagger)\mathbf{A}\mathbf{x}
\right\rangle
\right]
=
\mathbb{E}_{\mathbf{x}}
\!\left[
\left\langle
\mathbf{A}^\top \mathbf{B}\mathbb{E}[\mathbf{n}\mid\mathbf{x}],
(\mathbf{A}^\top\mathbf{B}-\mathbf{A}^\dagger)\mathbf{A}\mathbf{x}
\right\rangle
\right]
=0 .
\]

Consequently, under these assumptions one obtains the identity
\[
\mathbb{E}\|\mathbf{A}^\top \mathbf{B}\mathbf{y}-\mathbf{x}\|_2^2
=
\mathbb{E}\|\mathbf{A}^\top \mathbf{B}(\mathbf{y}-\mathbf{A}\mathbf{x})\|_2^2
+
\mathbb{E}\|(\mathbf{A}^\top\mathbf{B}-\mathbf{A}^\dagger)\mathbf{A}\mathbf{x}\|_2^2
+
\mathbb{E}\|(\mathbf{A}^\dagger\mathbf{A}-\mathbf{I})\mathbf{x}\|_2^2 .
\]

If the noise is not assumed to be centered or independent of $\mathbf{x}$,
the mixed term above cannot be guaranteed to vanish; in that case the same
derivation yields the inequality
\[
\mathcal{E}(\mathbf{B})
\;\le\;
\mathbb{E}\|\mathbf{A}^\top \mathbf{B}(\mathbf{y}-\mathbf{A}\mathbf{x})\|_2^2
+
\mathbb{E}\|(\mathbf{A}^\top\mathbf{B}-\mathbf{A}^\dagger)\mathbf{A}\mathbf{x}\|_2^2
+
\mathbb{E}\|(\mathbf{A}^\dagger\mathbf{A}-\mathbf{I})\mathbf{x}\|_2^2 ,
\]
which is the stated bound.
\end{proof}

\clearpage
\bibliographystyle{unsrt}   
\bibliography{sample}        

@book{Kak,
author = {Kak, Avinash C. and Slaney, Malcolm},
title = {Principles of Computerized Tomographic Imaging},
publisher = {Society for Industrial and Applied Mathematics},
year = {2001},
}

@book{Hensen,
title = "Computed Tomography: Algorithms, Insight, and Just Enough Theory",
editor = {Hansen, Per Christian and J{\o}rgensen, Jakob Sauer and Lionheart, William R. B.},
year = "2021",
language = "English",
isbn = "978-1-611976-66-3",
publisher = "Society for Industrial and Applied Mathematics"
}

@book{Gabor,
author = {Herman, Gabor T.},
title = {Fundamentals of Computerized Tomography: Image Reconstruction from Projections},
year = {2009},
isbn = {185233617X},
publisher = {Springer Publishing Company, Incorporated},
edition = {2nd}
}

@inproceedings{kudo,
  title={A very fast iterative algorithm for TV-regularized image reconstruction with applications to low-dose and few-view CT},
  author={Kudo, Hiroyuki and Yamazaki, Fukashi and Nemoto, Takuya and Takaki, Keita},
  booktitle={Developments in X-Ray Tomography X},
  volume={9967},
  pages={31--46},
  year={2016},
  organization={SPIE}
}

@article{bilgic,
  title={Fast image reconstruction with L2-regularization},
  author={Bilgic, Berkin and Chatnuntawech, Itthi and Fan, Audrey P and Setsompop, Kawin and Cauley, Stephen F and Wald, Lawrence L and Adalsteinsson, Elfar},
  journal={Journal of magnetic resonance imaging},
  volume={40},
  number={1},
  pages={181--191},
  year={2014},
  publisher={Wiley Online Library}
}

@article{geyer,
  title={State of the art: iterative CT reconstruction techniques},
  author={Geyer, Lucas L and Schoepf, U Joseph and Meinel, Felix G and Nance Jr, John W and Bastarrika, Gorka and Leipsic, Jonathon A and Paul, Narinder S and Rengo, Marco and Laghi, Andrea and De Cecco, Carlo N},
  journal={Radiology},
  volume={276},
  number={2},
  pages={339--357},
  year={2015},
  publisher={Radiological Society of North America}
}

@article{ART,
  title={Algebraic reconstruction techniques (ART) for three-dimensional electron microscopy and X-ray photography},
  author={Gordon, Richard and Bender, Robert and Herman, Gabor T},
  journal={Journal of theoretical Biology},
  volume={29},
  number={3},
  pages={471--481},
  year={1970},
  publisher={Elsevier}
}

@article{SIRT,
title = {Iterative methods for the three-dimensional reconstruction of an object from projections},
journal = {Journal of Theoretical Biology},
volume = {36},
number = {1},
pages = {105-117},
year = {1972},
issn = {0022-5193},
doi = {https://doi.org/10.1016/0022-5193(72)90180-4},
url = {https://www.sciencedirect.com/science/article/pii/0022519372901804},
author = {Peter Gilbert}
}

@article{Pan,
  title={Why do commercial CT scanners still employ traditional, filtered back-projection for image reconstruction?},
  author={Pan, Xiaochuan and Sidky, Emil Y and Vannier, Michael},
  journal={Inverse problems},
  volume={25},
  number={12},
  pages={123009},
  year={2009},
  publisher={IOP Publishing}
}

@article{FDK,
  title={Practical cone-beam algorithm},
  author={Feldkamp, Lee A and Davis, Lloyd C and Kress, James W},
  journal={Journal of the Optical Society of America A},
  volume={1},
  number={6},
  pages={612--619},
  year={1984},
  publisher={Optical Society of America}
}

@article{Sidky_2008,
year = {2008},
month = {aug},
volume = {53},
number = {17},
pages = {4777},
author = {Sidky, Emil Y and Pan, Xiaochuan},
title = {Image reconstruction in circular cone-beam computed tomography by constrained, total-variation minimization},
journal = {Physics in Medicine \& Biology}
}

@book{Buzug,
title = "Computed tomography: From photon statistics to modern cone-beam CT",
author = "Thorsten Buzug",
year = "2014",
month = jan,
day = "16",
doi = "10.1007/978-3-540-39408-2",
language = "English",
isbn = "978-3-540-39407-5",
publisher = "Springer Berlin Heidelberg"
}

@ARTICLE{Myagotin,
  author={Myagotin, Anton and Voropaev, Alexey and Helfen, Lukas and Hänschke, Daniel and Baumbach, Tilo},
  journal={IEEE Transactions on Image Processing}, 
  title={Efficient Volume Reconstruction for Parallel-Beam Computed Laminography by Filtered Backprojection on Multi-Core Clusters}, 
  year={2013},
  volume={22},
  number={12},
  pages={5348-5361},
  doi={10.1109/TIP.2013.2285600}
}

@inproceedings{Lauritsch,
  title={Theoretical framework for filtered back projection in tomosynthesis},
  author={Lauritsch, G{\"u}nter and H{\"a}rer, Wolfgang H},
  booktitle={Medical Imaging 1998: Image Processing},
  volume={3338},
  pages={1127--1137},
  year={1998},
  organization={SPIE}
}

@article{Yang,
  title={New reconstruction method for x-ray testing of multilayer printed circuit board},
  author={Yang, Min and Wang, Gao and Liu, Yongzhan},
  journal={Optical Engineering},
  volume={49},
  number={5},
  pages={056501--056501},
  year={2010},
  publisher={Society of Photo-Optical Instrumentation Engineers}
}

@article{Brien,
  title={Recent advances in X-ray cone-beam computed laminography},
  author={O’Brien, Neil S and Boardman, Richard P and Sinclair, Ian and Blumensath, Thomas},
  journal={Journal of X-ray Science and Technology},
  volume={24},
  number={5},
  pages={691--707},
  year={2016},
  publisher={IOS Press}
}

@book{Nesterov,
author = {Nesterov, Yurii},
title = {Introductory Lectures on Convex Optimization: A Basic Course},
year = {2014},
isbn = {1461346916},
publisher = {Springer Publishing Company, Incorporated},
edition = {1}
}

@article{Older,
doi = {10.1088/0031-9155/38/8/004},
url = {https://doi.org/10.1088/0031-9155/38/8/004},
year = {1993},
month = {aug},
volume = {38},
number = {8},
pages = {1051},
author = {J K Older and P C Johns},
title = {Matrix formulation of computed tomogram reconstruction},
journal = {Physics in Medicine \& Biology}
}

@article{Clack_1992,
doi = {10.1088/0031-9155/37/3/011},
url = {https://doi.org/10.1088/0031-9155/37/3/011},
year = {1992},
month = {mar},
volume = {37},
number = {3},
pages = {645},
author = {R Clack},
title = {Towards a complete description of three-dimensional filtered backprojection},
journal = {Physics in Medicine \& Biology}
}

@article{Zeng,
  title={A filtered backprojection algorithm with characteristics of the iterative Landweber algorithm},
  author={L. Zeng, Gengsheng},
  journal={Medical physics},
  volume={39},
  number={2},
  pages={603--607},
  year={2012},
  publisher={Wiley Online Library}
}

@article{Nielsen_2012,
doi = {10.1088/0031-9155/57/12/3915},
url = {https://doi.org/10.1088/0031-9155/57/12/3915},
year = {2012},
month = {may},
publisher = {IOP Publishing},
volume = {57},
number = {12},
pages = {3915},
author = {Nielsen, Tim and Hitziger, Sebastian and Grass, Michael and Iske, Armin},
title = {Filter calculation for x-ray tomosynthesis reconstruction},
journal = {Physics in Medicine \& Biology}
}

@article{Dan,
  author={Pelt, Daniël M. and Batenburg, Kees Joost},
  journal={IEEE Transactions on Image Processing}, 
  title={Improving Filtered Backprojection Reconstruction by Data-Dependent Filtering}, 
  year={2014},
  volume={23},
  number={11},
  pages={4750-4762},
  doi={10.1109/TIP.2014.2341971}}

@article{Joost,
  author={Batenburg, Kees Joost and Plantagie, Linda},
  journal={IEEE Transactions on Image Processing}, 
  title={Fast Approximation of Algebraic Reconstruction Methods for Tomography}, 
  year={2012},
  volume={21},
  number={8},
  pages={3648-3658},
  doi={10.1109/TIP.2012.2197012}}

@article{Dan_2,
  title={Improved tomographic reconstruction of large-scale real-world data by filter optimization},
  author={Pelt, Dani{\"e}l M and De Andrade, Vincent},
  journal={Advanced structural and chemical imaging},
  volume={2},
  number={1},
  pages={17},
  year={2016},
  publisher={Springer}
}

@article{blumensath,
  title={Backprojection inverse filtration for laminographic reconstruction},
  author={Blumensath, Thomas},
  journal={IET Image Processing},
  volume={12},
  number={9},
  pages={1541--1549},
  year={2018},
  publisher={Wiley Online Library}
}

@book{Saad,
author = {Saad, Yousef},
title = {Iterative Methods for Sparse Linear Systems},
publisher = {Society for Industrial and Applied Mathematics},
year = {2003},
doi = {10.1137/1.9780898718003},
edition   = {Second},
URL = {https://epubs.siam.org/doi/abs/10.1137/1.9780898718003},
eprint = {https://epubs.siam.org/doi/pdf/10.1137/1.9780898718003}
}

@article{hendriksen-2021-tomos,
author = {Hendriksen, Allard and Schut, Dirk and Palenstijn, Willem Jan and Viganò, Nicola and Kim, Jisoo and Pelt, Danië and van Leeuwen, Tristan and Batenburg, K. Joost},
title = {Tomosipo: Fast, Flexible, and Convenient {3D} Tomography for Complex Scanning Geometries in {Python}},
journal         = {Optics Express},
  year            = 2021,
  doi             = {10.1364/oe.439909},
  url             = {https://doi.org/10.1364/oe.439909},
  issn            = {1094-4087},
  month           = {Oct},
  publisher       = {The Optical Society},
}

@article{Astra,
author = {Wim van Aarle and Willem Jan Palenstijn and Jeroen Cant and Eline Janssens and Folkert Bleichrodt and Andrei Dabravolski and Jan De Beenhouwer and K. Joost Batenburg and Jan Sijbers},
journal = {Opt. Express},
number = {22},
pages = {25129--25147},
publisher = {Optica Publishing Group},
title = {Fast and flexible X-ray tomography using the ASTRA toolbox},
volume = {24},
month = {Oct},
year = {2016},
}

@article{Shepp_Logan_fil,
  author={Shepp, L. A. and Logan, B. F.},
  journal={IEEE Transactions on Nuclear Science}, 
  title={The Fourier reconstruction of a head section}, 
  year={1974},
  volume={21},
  number={3},
  pages={21-43},
  doi={10.1109/TNS.1974.6499235}}

@article{Kingma2014AdamAM,
  title={Adam: A Method for Stochastic Optimization},
  author={Diederik P. Kingma and Jimmy Ba},
  journal={CoRR},
  year={2014},
  volume={abs/1412.6980},
  url={https://api.semanticscholar.org/CorpusID:6628106}
}

@book{palamodov2016reconstruction,
  title={Reconstruction from Integral Data},
  author={Palamodov, V.},
  isbn={9781040219355},
  series={Chapman \& Hall/CRC Monographs and Research Notes in Mathematics},
  url={https://books.google.nl/books?id=tbCNEQAAQBAJ},
  year={2016},
  publisher={CRC Press}
}

@article{Bolzano,
title = {Bolzano and uniform continuity},
journal = {Historia Mathematica},
volume = {32},
number = {3},
pages = {303-311},
year = {2005},
issn = {0315-0860},
doi = {https://doi.org/10.1016/j.hm.2004.11.003},
url = {https://www.sciencedirect.com/science/article/pii/S0315086004000849},
author = {Paul Rusnock and Angus Kerr-Lawson}
}

@inproceedings{Petrica2024ANO,
  title={A note on identifiability for inverse problem based on observations},
  author={Marian Petrica and Ionel Popescu},
  year={2024},
  url={https://api.semanticscholar.org/CorpusID:271963011}
}

@book{nesterov2013introductory,
  title={Introductory lectures on convex optimization: A basic course},
  author={Nesterov, Yurii},
  volume={87},
  year={2013},
  publisher={Springer Science \& Business Media}
}

@book{villani2008optimal,
  title={Optimal transport: old and new},
  author={Villani, C{\'e}dric and others},
  volume={338},
  year={2008},
  publisher={Springer}
}

@InProceedings{Jin_2018_ECCV,
author = {Jin, Meiguang and Roth, Stefan and Favaro, Paolo},
title = {Normalized Blind Deconvolution},
booktitle = {Proceedings of the European Conference on Computer Vision (ECCV)},
month = {September},
year = {2018}
}

@article{Tan_DeepFBP,
  author={Tan, Xi and Liu, Xuan and Xiang, Kai and Wang, Jing and Tan, Shan},
  journal={IEEE Access}, 
  title={Deep Filtered Back Projection for CT Reconstruction}, 
  year={2024},
  volume={12},
  pages={20962-20972},
  doi={10.1109/ACCESS.2024.3357355}}

@inproceedings{faucris,
 author = {Sun, Yipeng and Schneider, Linda-Sophie and Fan, Fuxin and Thies, Mareike and Gu, Mingxuan and Mei, Siyuan and Zhou, Yuzhong and Bayer, Siming and Maier, Andreas},
 booktitle = {CT-Meeting},
 date = {2024-08-05/2024-08-09},
 doi = {10.48550/arXiv.2401.16039},
 faupublication = {yes},
 peerreviewed = {unknown},
 title = {{Data}-{Driven} {Filter} {Design} in {FBP}: {Transforming} {CT} {Reconstruction} with {Trainable} {Fourier} {Series}},
 venue = {Bamberg},
 year = {2024}
}
\end{document}